\documentclass[10pt]{article}

\usepackage[utf8]{inputenc}
\usepackage[T1]{fontenc}
\usepackage{fullpage}
\usepackage{enumitem}
\usepackage[shortcuts]{extdash}	
\usepackage{appendix}

\usepackage{xcolor}
\definecolor{mat}{HTML}{ffd6ad}

\usepackage{color}
\usepackage{bbm}					
\usepackage{relsize}				
\usepackage{amsmath}
\usepackage{amsthm}
\usepackage{amssymb}
\usepackage{amsfonts}
\usepackage{amscd}
\usepackage{stmaryrd}
\usepackage{MnSymbol}
\usepackage[mathcal]{euscript}			
\usepackage{tikz-cd}
\usepackage{hhline}
\usepackage{multirow}
\usepackage{tabularx}
\usepackage{hyperref}
\usepackage{xcolor}
\hypersetup{
    colorlinks=true,
    linkcolor=[RGB]{13 71 161},		
    citecolor=[RGB]{13 71 161},		
    urlcolor=[RGB]{13 71 161}		
}
\usepackage[capitalize,nameinlink]{cleveref}
\crefname{property}{Property}{Properties}
\Crefname{property}{Property}{Properties}
\crefname{fact}{Fact}{Facts}
\Crefname{fact}{Fact}{Facts}

\newtheorem{thm}{Theorem}[subsection]
\newtheorem{theorem}[thm]{Theorem}
\newtheorem{proposition}[thm]{Proposition}
\newtheorem{corollary}[thm]{Corollary}
\newtheorem{lemma}[thm]{Lemma}

\theoremstyle{definition}
\newtheorem{definition}[thm]{Definition}

\newtheorem{example}[thm]{Example}
\newtheorem{remark}[thm]{Remark}

\newtheorem{thmintro}{}
\newenvironment{thm-intro}[1]
  {\renewcommand\thethmintro{#1}\thmintro}
  {\endthmintro}

\newcommand{\ie}{{i.e. }}
\newcommand{\eg}{{e.g. }}
\newcommand{\oo}{$\infty$\=/}

\newcommand{\cA}{\mathcal{A}}
\newcommand{\cB}{\mathcal{B}}
\newcommand{\cC}{\mathcal{C}}
\newcommand{\cD}{\mathcal{D}}
\newcommand{\cE}{\mathcal{E}}
\newcommand{\cF}{\mathcal{F}}

\newcommand{\cM}{\mathcal{M}}

\newcommand{\cS}{\mathcal{S}}

\newcommand{\cU}{\mathcal{U}}

\newcommand{\sfE}{\mathsf{E}}
\newcommand{\sfF}{\mathsf{F}}

\newcommand{\sfS}{\mathsf{S}}

\newcommand{\sfU}{\mathsf{U}}

\newcommand{\NN}{\mathbb{N}}
\newcommand{\PP}{\mathbb{P}}
\newcommand{\QQ}{\mathbb{Q}}
\newcommand{\RR}{\mathbb{R}}
\newcommand{\ZZ}{\mathbb{Z}}

\newcommand\set{\mathsf{set}}
\newcommand\Set{\mathrm{\cS et}}

\newcommand\Fin[1]{{#1}_\mathrm{fin}}

\newcommand{\Cat}{\mathrm{\mathcal{C} at}}

\newcommand\pif[1]{{#1}_\pi}
\newcommand\Coh[1]{{#1}_{\overline\pi}}

\newcommand\Sh[1]{\mathrm{\cS h}\left(#1\right)}

\newcommand{\slice}[1]{_{/#1}}
\newcommand{\coslice}[1]{_{#1/}}

\newcommand{\Map}[2]{{\textrm{Map}\!\left(#1,#2\right)}}
\newcommand{\relMap}[3]{{\textrm{Hom}_{#1}\!\left(#2,#3\right)}}
\newcommand{\fun}[2]{\left[#1,#2\right]}
\newcommand{\Iso}[1]{{\textrm{Iso}\!\left(#1\right)}}

\newcommand{\pbh}[2]{\left\langle #1, #2 \right\rangle}
\newcommand{\Arr}[1]{{#1}^\rightarrow}
\newcommand\fperp{\upModels}

\DeclareMathOperator*{\colim}{colim}

\newcommand{\relP}[2]{{\mathcal{P}_{#1}\!\left(#2\right)}}

\newcommand{\op}{^{op}}

\renewcommand\index{_\bullet}
\newcommand\indexplus{_{\bullet+1}}
\newcommand\splx{^{\Delta\op}}
\newcommand\cart{_\mathsf{cart}}

\renewcommand\int{^\mathrm{int}}

\newcommand{\truncated}[1]{^{\leq #1}}

\newcommand\loc{^\mathsf{loc}}
\newcommand\cloc{^\mathsf{cloc}}

\newcommand{\pbmark}{\ar[dr, phantom, "\ulcorner" very near start, shift right=1ex]}
\newcommand{\pomark}{\ar[ul, phantom, "\lrcorner" very near start, shift right=1ex]}

\newcommand{\ot}{\leftarrow}
\newcommand{\tto}{{\begin{tikzcd}[ampersand replacement=\&]{}\ar[r]\&{}\end{tikzcd}}}
\newcommand{\stto}{{\begin{tikzcd}[ampersand replacement=\&, sep=small]{}\ar[r]\&{}\end{tikzcd}}}
\newcommand{\mto}{{\begin{tikzcd}[ampersand replacement=\&]{}\ar[r,mapsto]\&{}\end{tikzcd}}}
\newcommand{\xto}[1]{\xrightarrow {#1}}
\newcommand{\xot}[1]{\xleftarrow {#1}}
\newcommand{\subto}{\hookrightarrow}
\newcommand\surj{\twoheadrightarrow}
\newcommand\jrus{\twoheadleftarrow}

\newcommand\Loc[1]{\mathsf{Loc}(#1)}
\newcommand\BLoc[1]{\mathsf{BLoc}(#1)}
\newcommand\CLoc[1]{\mathsf{CLoc}(#1)}
\newcommand\BCLoc[1]{\mathsf{BCLoc}(#1)}

\begin{document}

\title{
The category of $\pi$\=/finite spaces%
\footnote{This material is based upon work supported by the Air Force Office of Scientific Research under award numbers FA9550-20-1-0305.}
}

\author{Mathieu Anel\footnote{Department of Philosophy, Carnegie Mellon University, mathieu.anel@protonmail.com}}

\maketitle

\begin{abstract}
We show that the category of truncated spaces with finite homotopy invariants
($\pi$\=/finite spaces) 
has many of the features expected of an elementary \oo topos.
It should be thought of as the natural higher analogue of the elementary 1-topos of finite sets, with which it shares several initiality properties.
The paper has also an appendix about univalent families in \oo pretopoi.
\end{abstract}

\setcounter{tocdepth}{2}
\tableofcontents

\section{Introduction}

Ever since the notion of \oo topos (which we shall call Grothendieck \oo topos in this discussion) started to be studied \cite{Simpson:topoi,Rezk:topos,TV:HAG1,Lurie:HTT}, the question of an elementary version of the notion has been around.
Such a notion would be to \oo topoi what Lawvere's elementary 1-topoi are to Grothendieck 1-topoi.
This question has become less academic with the discovery of the homotopical semantics of Martin-Löf's theory of dependent types and the introduction of the univalence axiom \cite{Awodey-Warren,Gambino-Garner,Voevodsky:lambda,Kapulkin-LeFanu-Lumsdaine,hottbook}.
The interpretation of logical types as homotopy types of spaces and identity types as path spaces has brought a deep and unexpected connection between logic and homotopy theory.
This connection has suggested a more precise content for the notion of elementary \oo topoi: they should be  the \oo categories which support an interpretation of dependent type theory with identity types and a univalent universe (aka homotopy type theory).

An axiomatization for elementary \oo topoi has been proposed in \cite{nlab:elementary} and an equivalent axiomatization has been developed in \cite{Rasekh:elementary}.
However examples are still scarce and somehow ad hoc \cite{Rasekh:elementary,Rasekh:FilterQuotients}. 
The purpose of this paper is to describe an example of an \oo category having many of the expected features required to interpret homotopy type theory.
Although this example does not verify all the axioms of the definition \cite{nlab:elementary,Rasekh:elementary}, it is nonetheless interesting because it is simple, concrete, and initial (or minimal) in several ways.

\begin{center}
*    
\end{center}

In 1-topos theory, the 1-category $\Set$ of (small) sets is a Grothendieck topos and the full subcategory $\Fin\Set$ of finite sets is an elementary topos.
When $\Set$ is generalized into the \oo category $\cS$ of spaces (\oo groupoids), the notion of finite sets has two natural generalizations:
\begin{enumerate}
\item {\it cell-finite spaces}, which are homotopy types of finite CW-complexes;
\item {\it $\pi$\=/finite spaces}, which are truncated homotopy types whose homotopy invariants are all finite (as a set or a group).
\end{enumerate}
The $\pi$\=/finite spaces are the {\it bounded coherent} objects in $\cS$ in the sense of \cite[Definition A.2.1.6 and Example A.2.1.7]{Lurie:SAG}.
They have been considered in relation to homotopy cardinality \cite{Baez-Dolan:cardinality,Baez:cardinality,Berman:cardinality,Yanovski:cardinality}, higher semiadditivity \cite{Lurie:ambidexterity,Harpaz:ambidexterity},
and pro-finite homotopy theory and higher Stone duality \cite[Appendix E]{Lurie:SAG}.

If $\cS$ is the \oo category of spaces, we denote by $\Fin\cS$ the full subcategory of cell-finite spaces and by $\pif\cS$ that of $\pi$\=/finite spaces.
These categories are essentially disjoint since their intersection is reduced to finite sets (\cref{prop:pi=omega})
\[
\begin{tikzcd}
\Fin\Set \ar[r, hook]\ar[d, hook] \pbmark & \pif\cS \ar[d, hook]\\
\Fin\cS \ar[r, hook] & \cS\,.
\end{tikzcd}
\]
Their stability properties are also very different (see \cref{figure:comparison}).
The \oo category $\Fin\cS$ is closed under finite colimits but not by fiber products.
Conversely, the \oo category $\pif\cS$ is closed under finite limits, finite sums, but not by pushouts.
In particular, the spheres $S^n$ ($n>0$) are cell-finite but not $\pi$\=/finite.
These properties make $\Fin\cS$ into a rather awkward object from the point of view of topos theory, where fiber products are fundamental.
The purpose of this paper is to show that $\pif\cS$, on the contrary, is very well behaved and definitely has some topos-like structure.

\medskip
We shall prove the following properties of $\pif\cS$:
\begin{enumerate}[label=(\arabic*)]
\item \label[property]{prop-coh:1} it is a lex \oo category (\cref{prop:coh-lex});
\item \label[property]{prop-coh:2} which is extensive (\ie finite sums exist, and are universal and disjoint, see \cref{prop:coh-extensive});
\item \label[property]{prop-coh:3} and exact (\ie quotients of Segal groupoids objects exist, and are universal and effective, see \cref{prop:Segal-gpd});
\item \label[property]{prop-coh:4} $\pif\cS$ has all truncation (Postnikov) modalities (\cref{prop:coh-truncation});
\item \label[property]{prop-coh:5} it is locally cartesian closed (\cref{thm:CohLCC});
\item \label[property]{prop-coh:8} its universe $\pif U$ (which lives in $\cS$) is a countable sum of $\pi$\=/finite spaces (\cref{thm:universe});
\item \label[property]{prop-coh:6} $\pif\cS$ has enough univalent maps (\cref{thm:univalent}), and they can be chosen closed under diagonals, aka identity types (\cref{thm:univalent:id});
\item \label[property]{prop-coh:7} it has a subobject classifier, which is Boolean (\cref{prop:CohOmega}).
\end{enumerate}

\Cref{prop-coh:1,prop-coh:2,prop-coh:3} make $\pif\cS$ into an \oo pretopos in the sense of \cite[Appendix A]{Lurie:SAG} where it is mentioned as an example.
In fact, we shall see that
\begin{enumerate}[label=(\arabic*)]
\setcounter{enumi}{8}
\item \label[property]{prop-coh:9} $\pif\cS$ is the initial \oo pretopos (\cref{thm:initial}).
\end{enumerate}
We also show a couple of stronger universal properties:
\begin{enumerate}[label=(\arabic*)]
\setcounter{enumi}{9}
\item \label[property]{prop-coh:10bis} $\pif\cS$ is the initial ``locally cartesian closed \oo pretopos'' (see \cref{thm:initial-LCCpretopos} for a precise statement);
\item \label[property]{prop-coh:10} $\pif\cS$ is the initial ``locally cartesian closed \oo pretopos with a Boolean subobject classifier'' (see \cref{thm:initial-boolean} for a precise statement).
\end{enumerate}
\Cref{prop-coh:4} is a consequence of the \oo pretopos structure, but it is obvious here.
The proof of that $\pif\cS$ is locally cartesian closed (\cref{prop-coh:5}) will come after a study of the descent properties of $\pif\cS$, which form the bulk of \cref{sec:LCC}. 
The main tool will be the folklore result characterizing $\pi$\=/finite spaces as realization of Kan complexes with values in finite sets (\cref{prop:resolution}).
From there, the construction of the universe and the study of univalent families (\cref{prop-coh:8,prop-coh:6}) are fairly straightforward.

\medskip
Altogether, this provides the \oo category $\pif\cS$ with almost all of the properties of the notion of elementary \oo topos of \cite{nlab:elementary,Rasekh:elementary}, but
\begin{enumerate}[label=(\alph*)]
\item \label[fact]{fact:1} $\pif\cS$ does not have all pushouts (\eg the spheres $S^n$ for $n>0$ are not $\pi$\=/finite, see \cref{prop:pushout}),
\item \label[fact]{fact:2} and it does not have a hierarchy of univalent families closed under dependent sums and/or dependent products (\cref{thm:no-universe}).
\end{enumerate}

In an elementary 1-topos, the existence of pushouts can be deduced from the existence of finite limits, exponentials and the subobject classifier.
\Cref{fact:1} shows that $\pif\cS$, viewed as a universe in $\cS$, is a counter-example in higher categories, providing a negative answer to a question of Awodey.
Nonetheless, it is established in \cite{Frey-Rasekh}, that finite sums can be build from finite limits, exponentials and the subobject classifier. 
It seems reasonable that the pushouts where one leg is a monomorphism can also be produced this way.
\Cref{fact:2} is essentially due to the fact that there are no inaccessible cardinals between 2 and $\omega$. 
This is related to the minimality properties of $\pif\cS$.

\medskip
In connection to homotopy type theory, both \cref{fact:1,fact:2} are considered as no-go, and, since the first version of this note, the consensus seems to be that $\pif\cS$ should not be an elementary topos.
It remains that $\pif\cS$ does provides a non-trivial univalent family $\pif U'\to \pif U$ in $\cS$ which is closed under diagonals, dependent sums, dependent products, finite sums, finite products, and quotients of groupoidal equivalence relations (\cref{thm:universe}).
This universe also exists in the subcategory of $\kappa$-small spaces for every inaccessible cardinal $\kappa>\omega$.
The universe $\pif U$ is not the smallest universe of $\cS$ closed under dependent sums and products (since the subobject classifier $\Omega=2$ or the 1-type $\coprod_n B\mathfrak S_n$ of finite sets are also examples), but \cref{prop-coh:10bis} can be reformulated by saying that $\pif U$ is the minimal universe of $\cS$ containing the universe $\coprod_n B\mathfrak S_n$ of finite sets, and closed under quotients of groupoidal equivalence relations (see \cref{cor:initial-LCCpretopos}).

\paragraph{Acknowledgments}
I thank 
Carlo Angiuli, 
Steve Awodey, 
Reid Barton, 
Jonas Frey, 
André Joyal, 
Nima Rasekh, 
Mike Shulman, and 
Andrew Swan 
for many discussions about 1-pretopoi, elementary 1-topoi, type theory, and their comments about earlier drafts.
\Cref{thm:initial-boolean} was suggested by Steve Awodey.
I thank Tim Holzschuh and Maxime Ramzi for pointing out a mistake in the ``building blocks'' of $\pif\cS$.
I learned the theory of \oo pretopos in the Appendix A of Jacob Lurie's book \cite{Lurie:SAG}, many techniques and results are taken from there.
I thank also the anonymous referee for many remarks that have helped improve the paper.

\paragraph{Convention}
This paper is written in the language of \oo categories but we shall drop all ``$\infty$-'' prefixes and call higher categorical notions by their classical name (category always means \oo category, topos means \oo topos, colimit always means \oo colimit, pullback always means \oo pullback, etc.)
When $n$-categories and $n$-categorical notions will be required for $n<\infty$, we shall use an explicit ``$n$-'' prefix.
We refer to \cite{Lurie:HTT,Cisinski,Riehl-Verity:EICT} for basics on \oo category theory.
All the arguments of the paper are formulated in terms that make sense in any model of \oo category theory (limits, colimits, exactness properties, adjunctions...).

\section{\texorpdfstring{$\pi$}{pi}-finite spaces}

\subsection{Definition and characterizations}

We shall say that a space $X$ is {\it unbounded $\pi$\=/finite} if $\pi_0(X)$ is a finite set and all $\pi_n(X,x)$ ($n>0$) are finite groups, for all choices of base point.
We shall say that a space $X$ is {\it $\pi$\=/finite} if it is moreover $n$-truncated for some $n$.
Unbounded $\pi$\=/finite ($\pi$\=/finite) spaces are the coherent (bounded coherent) objects of the category $\cS$ \cite[Example A.2.1.7.]{Lurie:SAG}.
We denote by $\Coh\cS$ and $\pif\cS$ the subcategory of $\cS$ spanned by unbounded $\pi$\=/finite and $\pi$\=/finite spaces.

In the case of spaces, the coherence condition can be understood as a higher analogue of the notion of Kuratowski finite object.
We say that a space $X$ is {\it finitely covered} if there exist a map $E\to X$ where $E$ is a finite set and which is surjective on $\pi_0$.
We say that a map $X\to Y$ is {\it finitely covered} if all its fiber are finitely covered.
Recall that the diagonal of a map $f:X\to Y$ is the map $\Delta f:X\to X\times_YX$.
The higher diagonals are defined by $\Delta^{n+1}f = \Delta(\Delta^n f)$.
When $Y=1$ is the point, we have $\Delta^{n+1} X :=\Delta^{n+1} (X\to 1) = X \to X^{S^n}$.

\begin{proposition}[Kuratowski characterization]
A space $X$ is unbounded $\pi$\=/finite if and only if all its diagonals $\Delta^{n+1}X$ are finitely covered.
\end{proposition}
\begin{proof}
For any base point $x$ in $X$, the sets $\pi_n(X,x)$ is finite if and only if $\Omega^n_xX$ is finitely covered.
The result follows from the fact that the fibers of $\Delta^{n+1}X:X \to X^{S^n}$ are exactly the loop spaces $\Omega^{n+1}X$.
\end{proof}

\medskip
\noindent Examples of $\pi$\=/finite spaces:
\begin{itemize}
\item any finite set (including the empty set $\emptyset$ and singleton 1);
\item $\RR \PP^\infty = B\ZZ_2$ (= universe of sets of cardinal 2);
\item $\coprod_{k\leq n} B\mathfrak S_k$ (= universe of sets or cardinal $\leq n$);
\item the classifying space $BG$ of a finite group $G$;
\item more generally, the realization of any finite groupoid ($G_1\rightrightarrows G_0$ in $\Fin\Set$);
\item the Eilenberg--Mac~Lane space $K(G,n)$ for a finite group $G$.
\item the ``local coefficient'' Eilenberg--Mac~Lane spaces $K(G,A,n)$ where $G$ is a finite group, $A$ a finite $G$-module (which can be defined as the quotient of $K(A,n)$ by the action of $G$) \cite[VI.3]{Goerss-Jardine}.

\end{itemize}

\noindent Examples of unbounded $\pi$\=/finite spaces:
\begin{itemize}
\item By Serre finiteness theorem, sufficiently large loop spaces of the sphere $S^n$ are spaces with finite homotopy groups ($\Omega^{m}S^{n}$ for $m>n$ when $n$ is odd, and for $m\geq 2n$ when $n$ is even);
\item The realization of any Kan complex with values in finite sets (see \cref{prop:resolution}).
\end{itemize}


\subsection{Elementary properties}

This section proves \cref{prop-coh:1,prop-coh:2,prop-coh:4} of $\pif\cS$.

\begin{lemma}
\label{lem:coh-subobject}
Any subspace of a $\pi$\=/finite space is $\pi$\=/finite.
\end{lemma}
\begin{proof}
A subspace is determined by a subset of connected components.
Hence, the $\pi_0$ is finite and so are the higher homotopy groups.
\end{proof}

\begin{lemma}
\label{lem:coh-extensive}
The category $\pif\cS$ has finite sums and the inclusion $\pif\cS\subset \cS$ preserves them.
\end{lemma}
\begin{proof}
The initial object of $\cS$ is $\pi$\=/finite.
Let $X$ and $Y$ be two $\pi$\=/finite spaces, then the sum $X+Y$ (computed in $\cS$) is $\pi$\=/finite and provide a sum for $X$ and $Y$ in $\pif\cS$.
\end{proof}

The following result proves that the category $\pif\cS$ is closed under fibers, extensions, and quotients (see \cref{prop:Segal}).
\begin{proposition}
\label{prop:342}
\label{lem:coh-lex}
\label{lem:coh-sigma}
Consider a cartesian square
\[
\begin{tikzcd}
Z \ar[r]\ar[d] \pbmark & X \ar[d]\\
1 \ar[r] & Y
\end{tikzcd}
\]
where $Y$ is a connected space.
Then, if any two of $X$, $Y$, or $Z$ are $\pi$\=/finite, so is the third.
\end{proposition}
\begin{proof}
By working componentwise, we can assume that $X$ is connected.
We choose an arbitrary base point $z$ in $Z$, we denote $x$ its image in $X$ and $y$ its image in $Y$.
We consider the long exact sequence of homotopy invariants:
\[
\dots \pi_2(Z,z) \to \pi_2(X,x) \to \pi_2(Y,y) \to 
\pi_1(Z,z) \to \pi_1(X,x) \to \pi_1(Y,y) \to 
\pi_0(Z) \to 1\,.
\]
We prove the result in case where $X$ and $Y$ are assumed in $\pif\cS$.
The map $\pi_1(Y,y) \to \pi_0(Z,z)$ is surjective, this proves that $\pi_0(Z,z)$ is finite.
For $n>0$, we get a short exact sequence $K\to \pi_n(Z,z)\to Q$
where $K$ is the kernel of $\pi_n(Z,z)\to \pi_n(X,x)$,
and $Q$ is the quotient of the map $\pi_{n+1}(Y,y)\to \pi_n(Z,z)$.
$K$ is a subgroup of a finite group, $Q$ is a quotient of a finite group, hence they are both finite.
Then $\pi_n(Z,z)$ is finite since, as a set, it is in bijection with $K\times Q$.
Since the base point of $Z$ was arbitrary, this proves that $Z$ is in $\pif\cS$.
The argument is similar in the two other cases.
\end{proof}

\begin{proposition}[Finite limits]
\label{prop:coh-lex}
The category $\pif\cS$ has finite limits (in particular loop spaces) and they are preserved by the inclusion $\pif\cS\subset \cS$.
\end{proposition}
\begin{proof}
The point is $\pi$\=/finite.
The statement for binary products is direct from the formula $\pi_n(X\times Y,(x,y)) = \pi_n(X,x)\times\pi_n(Y,y)$.
We need only to check fiber products.
Given a diagram $X\to Y\ot Y'$ in $\pif\cS$, we want to prove that $X\times_YY'$ is in $\pif\cS$.
We can work over each connected component of $Y$ separately and assume that $Y$ is connected.
Using \cref{lem:coh-sigma}, it is enough to prove that the fibers of the map $X\times_YY'\to Y'$ are in $\pif\cS$.
But these fibers are fibers of the map $X\to Y$, which are in $\pif\cS$ by \cref{lem:coh-lex}.
\end{proof}

The $n$-truncation of a $\pi$\=/finite space is clearly $\pi$\=/finite.
The next result proves that the results holds not only for objects but also maps.

\begin{proposition}[Postnikov truncations]
\label{prop:coh-truncation}
Let $X\to Y$ be a morphism in $\pif\cS$, and $X\to Z\to Y$ its factorization (computed in $\cS$) into an $n$-connected maps followed by an $n$-truncated map.
Then, the space $Z$ is $\pi$\=/finite.
\end{proposition}
\begin{proof}
We can work over each connected component of $Y$ separately and assume that $Y$ is connected.
The fiber of $Z\to Y$ is the $n$-truncation of the fiber of $X\to Y$,
and the result follows from \Cref{lem:coh-sigma}.
\end{proof}

\begin{definition}[Descent]
\label{def:descent}    
Let $\cC$ be a category with finite limits and all colimits indexed by some small category $I$.
For any diagram $X:I\to \cC$, the base change functor induces a functor $P:\cC\slice {\colim X_i} \to \lim_i\cC\slice {X_i}$
and the colimit functor induces a left adjoint $C:\lim_i\cC\slice {X_i} \to \cC\slice {\colim X_i}$.
We say that an $I$-colimit is {\it universal} (resp. {\it effective}) if $C$ (resp. $P$) is a fully faithful functor.
We say that an $I$-colimit has {\it descent} if it is universal and effective \cite[6.1.3, 6.1.8]{Lurie:HTT}.
\end{definition}

If $I$ is a set, effectivity corresponds to the disjunction of sums, and descent to their extensivity.
All colimits have descent in $\cS$ \cite[Theorem~6.1.3.9 \& Lemma~6.1.3.14 ]{Lurie:HTT}.

\begin{lemma}[Descent]
\label{lem:descent}
If $\cC\subset\cS$ is a subcategory closed under finite limits such that the colimit of some small diagram $X:I\to \cC$ exist in $\cC$ and is preserved by $\cC\subset\cS$, then this colimit has descent in $\cC$.
\end{lemma}
\begin{proof}
Let $X:I\to \cC$ be such a diagram.
By assumption, the adjunction $\cC\slice {\colim X_i} \rightleftarrows \lim\,\cC\slice {X_i}$
is the restriction to $\cC$ of the adjunction $\cS\slice {\colim X_i} \rightleftarrows \lim\,\cS\slice {X_i}$ which is an equivalence of categories by descent in $\cS$.
Hence so is the adjunction $\cC\slice {\colim X_i} \rightleftarrows \lim\,\cC\slice {X_i}$.
\end{proof}

\begin{proposition}[Extensivity]
\label{prop:coh-extensive}
Finite sums in $\pif\cS$ have descent (are disjoint and universal).
\end{proposition}
\begin{proof}
We saw that the inclusion $\pif\cS\subset \cS$ preserves finite limits and finite sums in \cref{prop:coh-lex} and \cref{lem:coh-extensive}.
The result follows from \Cref{lem:descent}.
\end{proof}

\begin{proposition}
\label{prop:retract}
The category $\pif\cS$ is idempotent complete.
\end{proposition}
\begin{proof}
Let $Y$ be a retract of a $\pi$\=/finite space $X$, then, for every base point $y$,
$\pi_n(Y,y)$ is a retract of $\pi_n(X,y)$, hence finite and eventually null.
\end{proof}

\begin{remark}
\label{rem:coh}
All the previous results are also true for the category $\Coh\cS$ of unbounded $\pi$\=/finite spaces (with the same proofs).
\end{remark}

Using Postnikov towers, any unbounded connected $\pi$\=/finite space can be built by successive pullbacks of maps $K(G,1)\to K(G,A,n)$ where $G$ is a finite group and $A$ a finite $G$-module \cite[VI.3]{Goerss-Jardine}.
\Cref{prop:coh-truncation} ensures that all steps of the construction are in $\pif\cS$.
This gives the following result.

\begin{proposition}
\label{prop:coh-block}
The category of $\pi$\=/finite spaces is the smallest subcategory of $\cS$ closed under finite sums, finite limits, and containing all ``local coefficient'' Eilenberg--Mac~Lane spaces $K(G,A,n)$ for $G$ a finite group and $A$ a finite $G$-module.
\end{proposition}

\subsection{Absence of pushouts}

This section proves that $\pif\cS$ does not have all pushouts.
We know already that the inclusion $\pif\cS\subset \cS$ cannot preserve pushouts since $S^1 = 1\cup_{S^0}1$ is not $\pi$\=/finite.
But that does not prevent pushouts to exist in $\pif\cS$.
We will see that it is indeed the case by proving that the pushout $1\ot 2\to 1$ (where $2=1+1$), classically equal to $S^1$, does not exist in $\pif\cS$.

\begin{proposition}
\label{prop:pushout}
The pushout of the diagram $1\ot 2\to 1$ is not representable in $\pif\cS$.
\end{proposition}
\begin{proof}
This pushout is by definition the object representing the free loop space functor 
\begin{align*}
FL:\pif\cS &\tto \cS\\
X &\mto X\times_{X\times X}X = X^{S^1}
\end{align*}
Let $H$ be a $\pi$\=/finite space representing $FL$.
Then, for any unbounded $\pi$\=/finite space $X$, we have a natural equivalence $\Map H X = \Map {S^1} X$.
This is equivalent to the data of a map $S^1\to H$ presenting $H$ as the reflection of $S^1$ in the subcategory $\pif\cS\subset \cS$.
First we can deduce that $H$ has to be connected.
Indeed, for the discrete space $X=2$, we have
\[
\Map H 2 = 2^{\pi_0(H)}
\qquad \text{and} \qquad
\Map {S^1} 2 = 2
\]
hence $\pi_0(H)=1$.

Recall that for $(X,x)$ and $(Y,y)$ two pointed spaces, the space of pointed maps is defined by the fiber product
\[
\begin{tikzcd}
\relMap * {(Y,x)} {(X,x)} \ar[r]\ar[d] \pbmark & \Map  Y X \ar[d,"\Map y X"]\\
1 \ar[r,"x"] & X
\end{tikzcd}
\]
We fix a base point $s$ in $S^1$, and consider its image $h$ by the map $S^1\to H$.
The map $S^1\to H$ induces an equivalence 
\begin{equation}
\label{eqn:po}
\relMap * {(H,h)} {(X,x)}  = \relMap * {(S^1,s)} {(X,x)} = \Omega_x X
\end{equation}
The space $H$ being connected its 1-truncation is a space $BG$ for some finite group $G$.
We consider the additive group $\ZZ/p\ZZ$ for $p$ a prime number prime to the order of $G$.
Then, the only group morphism $G\to \ZZ/p\ZZ$ is the constant one.
We put $X = B\ZZ/p\ZZ$.
Using the equivalence between pointed connected 1-type and discrete groups, we get
\[
\relMap * {(H,h)} {(X,x)} = \relMap * {(BG,h)} {(B\ZZ/p\ZZ,x)} = \relMap {Gp} G {\ZZ/p\ZZ} = 1.
\]
But, on the other side, we have
\[
\Omega_x X = \ZZ/p\ZZ \not=1
\]
This contradicts \eqref{eqn:po} and shows that $H$ cannot exist.
\end{proof}

\begin{remark}
The argument can be adapted to show that $S^n$ ($n\geq 1$) does not admit a reflection $H$ into $\pif\cS$.
First, $H$ must be $(n-1)$-connected since $\relMap * H X = \Omega^nX = 1$ for every $(n-1)$-truncated space (and any choice of base point $h$ on $H$).
Then using $X=K(\ZZ/p\ZZ,n)$, where $p$ is prime to the cardinality of $G:=\pi_n(H,h)$, we get the contradiction $\ZZ/p\ZZ = \Omega^n X = \relMap * H X = \relMap * {K(G,n)} {K(\ZZ/p\ZZ,n)} = 1$.
It is likely that any nontrivial connected cell-finite space does not admits a reflection either.
\end{remark}

Not all pushouts exists in $\pif\cS$ but some do.
\begin{proposition}
The pushouts of spans where one of leg is a monomorphism exist in $\pif\cS$.
\end{proposition}
\begin{proof}
All monomorphisms are split in $\cS$, hence in $\pif\cS$: if $X\to Y$ is a monomorphism in $\pif\cS$ then $Y$ is isomorphic to $X\coprod Y'$ for some $Y'$ in $\pif\cS$.
Then, given a span $Y\ot X\to Z$, the pushout is $Z\coprod Y'$ which is in $\pif\cS$.
\end{proof}

\subsection{Comparison with cell-finite spaces}

We recall without proof some properties of cell-finite spaces to compare them with the $\pi$\=/finite ones.
The comparison between the two categories is summarized in \cref{figure:comparison}.

\medskip
Recall that we call a space {\it cell-finite} if it is the homotopy type of a finite CW-complex,
or, equivalently, the realization of a simplicial set with only a finite number of non-degenerate simplices.
More intrinsically, the category of cell-finite spaces can be defined as the smallest subcategory of $\cS$ containing $\emptyset$ and 1 (or the whole of $\Fin\Set$) and closed under pushouts.
(We shall see in \cref{prop:Segal-gpd}, a similar characterization of unbounded $\pi$\=/finite spaces.) 
All spheres $S^n$ ($n\geq -1$) are cell-finite, and any cell-finite space can be built with a 
finite chain of cell attachments
\[
\begin{tikzcd}[ampersand replacement=\&]
S^n \ar[d]\ar[r] \& X_n\ar[d]\\
1 \ar[r]\& X_{n+1}\pomark \, .
\end{tikzcd}
\]
This is to be contrasted with \cref{prop:coh-block}.
Any subspace of an cell-finite space is cell-finite.
Any finite sums or finite product of cell-finite spaces is cell-finite.
But $\Fin\cS\subset \cS$ is not closed under finite limits since $\Omega S^1=\ZZ$ is not cell-finite.
It is also not closed under retracts \cite[Remark 5.4.1.6]{Lurie:HTT}.
\Cref{figure:comparison} summarizes the comparison between cell-finite and $\pi$\=/finite spaces.
A funny fact is that, $\pif\cS$ being closed under finite limits, it is cotensored over $\Fin\cS$:

\begin{lemma}
\label{lem:cotensor}
The mapping space functor $\mathrm{Map}:\cS\op\times \cS\to \cS$ restricts into a functor $\Fin\cS\op\times \pif\cS\to \pif\cS$.
\end{lemma}
\begin{proof}
Let $K$ be an cell-finite space and $X$ be a $\pi$\=/finite space, then $X^K$ is a finite limit of copies of $X$, hence in $\pif\cS$ by \cref{prop:coh-lex}.
\end{proof}

\begin{proposition}
\label{prop:pi=omega}
A space is $\pi$\=/finite and cell-finite if and only if it is a finite set.
\end{proposition}
\begin{proof}
Clearly $\Fin\Set \subseteq \pif\cS\cap \Fin\cS$.
Any space in $\pif\cS\cap \Fin\cS$ has a finite number of connected components, to prove the expected equality, it is then enough to show that the spaces which are connected, $\pi$\=/finite and cell-finite are contractible.

We will use the theory of fiberwise orthogonality and acyclic classes of \cite[3.2]{ABFJ:HS}.
We say that two maps $u:A\to B$ and $f:X\to Y$ are orthogonal (denoted $u\perp f$) if the map $\pbh u f := X^B \to X^A \times_{Y^A} Y^B$ is invertible.
We say that two maps $u$ and $f$ are fiberwise orthogonal (denoted $u\fperp f$) if $u'\perp f$ for every base change $u':A'\to B'$ of $u$.
Notice that for two objects $X$ and $Y$, we have $(X\to 1) \fperp (Y\to 1)$ if and only if $(X\to 1) \perp (Y\to 1)$ if and only if $Y\to Y^X$ is invertible.
Miller's theorem (aka Sullivan's conjecture) \cite{Miller:Sullivan} says that the canonical map $Y\to \Map {BG} Y$ is invertible (we are considering unpointed maps) when $Y$ is a cell-finite space and $G$ a finite group.
This is equivalent to say that $(BG\to 1) \fperp (Y\to 1)$.

For any class of maps $\cB$, the class $^\fperp\cB:= \{u\ |\ \forall f\in \cB, u\fperp f\}$ is acyclic, which means that 
it contains all isomorphisms, it is closed under composition, base change, and under small colimits (computed in $\Arr\cS$).
By Miller's theorem, the acyclic class $\cM = ^\fperp \{X\to 1\ |\ X\in \Fin\cS\}$ contains all maps $BG\to 1$ where $G$ is a finite group.
Let us see that $\cM$ contains the maps $Z\to 1$ for $Z$ an arbitrary connected $\pi$\=/finite space.
We know that all $BG$ are in $\cM$ for all finite groups.
An induction shows that, for all finite abelian groups $G$ and for all $n\geq 1$, the spaces $K(G,n)$ are in $\cM$.
We've seen this is the case for $n=1$.
Suppose this is the case for all $n\leq m$, we show that this is also the case for $n= m+1$:
for any finite abelian group $G$, the space $K(G,m+1)= B(K(G,m))$ is the colimit of the simplicial diagram encoding the group structure of $K(G,m)$,
thus the maps $K(G,m+1)\to 1$ belongs to $\cM$ since acyclic classes are closed under colimits.
Therefore all maps whose fibers are $K(G,n)$ are in $\cM$.
By closure of $\cM$ under composition, all maps whose fibers are connected $\pi$\=/finite spaces are also in $\cM$.
In particular, any map $X\to 1$ where $X$ is a connected $\pi$\=/finite space is in $\cM$.
So, if $X$ is a space which is connected, $\pi$\=/finite and cell-finite, by orthogonality between $\cM$ and cell-finite spaces, the diagonal map $X\to X^X$ must be invertible.
Thus the identity of $X$ is homotopic to a constant map and $X$ is contractible.
\end{proof}

\begin{remark}[Modality generated by finite $\pi$-spaces]
Recall from \cite[3.2]{ABFJ:HS} that a modality is a (unique) factorization system stable under base change.
The proof of \cref{prop:pi=omega} shows that there exists a non-trivial modality $(\cA,\cB)$ on the category of $\cS$ where the left class $\cA$ contains all maps in $\cS$ whose fibers are connected $\pi$\=/finite spaces, and where the right class $\cB$ contains all maps in $\cS$ whose fibers are cell-finite spaces.
The corresponding factorization is unknown to us.
\end{remark}

\begin{remark}[Generalization to higher cardinals]
\label{sec:higher-cardinal}
The notions of cell-finite and $\pi$\=/finite spaces (and more generally that of compact and coherent object in a topos) rely implicitly on the notion of finite sets, that is $\omega$-small sets.
It can therefore be generalized by replacing $\omega$ with a non countable larger regular ordinal $\kappa$.
If we do so, then the notion of $\kappa$-small and $\kappa$\=/finite spaces do coincide.
Only for $\omega$ are the two notions different.
An explanation is the following: the completion of a simplicial set with a values in finite sets (a fortiori having a finite set of non-degenerate simplices) into a Kan complex has values in  countable sets (using the $Ex^\infty$ fibrant replacement, each iteration of $Ex$ functor stays with finite values, but we need a countable iteration).
When $\kappa>\omega$, the $Ex^\infty$ of a simplicial set with values in $\kappa$-small sets stay with values in $\kappa$-small sets.
\end{remark}

\begin{remark}[Closure of $\pif\cS$ for pushouts]
\label{sec:closure-pushout}
A natural question is to identify the smallest category $\cS_?$ containing $\pif\cS$ and $\Fin\Set$ which is closed under finite limits and finite colimits.
An upper bound is given by the subcategory $\cS_\mathrm{count.}$ of $\cS$ spanned by realizations of countable simplicial sets, since it is closed under finite limits, finite colimits and contains $\pif\cS\cup\Fin\Set$.
Countable sets will be in $\cS_?$, since it must contains $\ZZ=\Omega S^1$, but not the Hom sets between countable sets since they are sets outside $\cS_\mathrm{count.}$
This implies that $\cS_?$ cannot be cartesian closed (and thus locally cartesian closed).
This category has been studied by Berman \cite{Berman:cardinality}.
\end{remark}

\begin{table}[htbp]
\begin{center}
\caption{Comparison between cell-finite and $\pi$\=/finite spaces.}
\label{figure:comparison}

\medskip
\renewcommand{\arraystretch}{1.6}
\begin{tabularx}{.8\textwidth}{
|>{\hsize=1\hsize\linewidth=\hsize\centering\arraybackslash}X
|>{\hsize=.6\hsize\linewidth=\hsize\centering\arraybackslash}X
|>{\hsize=1.4\hsize\linewidth=\hsize\centering\arraybackslash}X|}
\cline{2-3}
\multicolumn{1}{c|}{}& $\Fin\cS$ & $\pif\cS$ \\
\hline
finite $+$ and $\times$ & yes & yes \\
\hline
subspaces & yes & yes \\
\hline
pushouts & yes & no \\
\hline
fiber products & no & yes  \\
\hline
loop spaces $\Omega$ & no & yes  \\
\hline
truncations & no & yes  \\
\hline
retracts & no & yes \\
\hline
building blocks & $S^n\to 1$ by~pushouts & $1\to K(G,A,n)$ by~fiber~products \\
\hline
compactness properties & finite spaces are compact in $\cS$ & $n$-truncated $\pi$\=/finite spaces are compact in $\cS\truncated n$ (but not in $\cS$) \\
\hline
Euler characteristic 
& in $\ZZ = \NN[-1]$ & in $\QQ_{>0} = \NN[{1\over 2},{1\over 3},\dots]$ 
\cite{Baez-Dolan:cardinality,Baez:cardinality,Berman:cardinality,Yanovski:cardinality}  \\
\hline
Other properties
&---&
higher semiadditivity
\cite{Lurie:ambidexterity,Harpaz:ambidexterity} \\
&&
cotensored over $\Fin\cS$ (\cref{lem:cotensor}) \\
\hline
\end{tabularx}
\end{center}
\end{table}

\section{Descent and cartesian closure}
\label{sec:LCC}
The main result of this section is that $\pif\cS$ is locally cartesian closed (\cref{thm:CohLCC}).
We start with some recollections about groupoid objects and descent.

\subsection{Simplicial spaces}

This section introduces definitions and constructions useful in the following sections.

\medskip
A map $f:X\to Y$ in $\cS$ is called {\it surjective} if the map $\pi_0(f):\pi_0(X)\to \pi_0(Y)$ is surjective in $\Set$.
An object $E$ of $\cS$ is called {\it projective} if for any surjective map $X\to Y$, the map $X^E \to Y^E$ is surjective.
We shall say that a map $f:X\to Y$ in $\cS$ is {\it projective} if it is of the type $X\to X+E$ where $E$ is a projective object.

\begin{lemma}[Enough projective]
\label{lem:choice}
In $\cS$, the projective objects are the sets (i.e. the 0-truncated spaces).
Moreover, $\cS$ has enough projective objects in the sense that any map $X\to Y$ can be factored into $X\to X+E\to Y$ into a projective map followed by a surjection.
The projective maps and surjections form a (non functorial) weak factorization system on $\cS$.
\end{lemma}
\begin{proof}
Clearly, any set is projective.
Conversely, surjective maps being closed under base change, $E$ is projective if and only if any surjective map $X \to E$ splits.
Using that any object $Y$ of $\cS$ admits (non functorially) a surjective map $X\to Y$ from a set
(in the model of $\cS$ with topological spaces, $E$ can be the set of points of $Y$; in the model with simplicial sets, $E$ can be the set of 0-simplices).
We see that a projective object $E$ must be a retract of a set, hence a set.
This proves the first statement.
For the factorization, it is enough to choose any projective cover $E\to Y$ from a set and consider $X\to X+E\to Y$.
The non functoriality is due to that of the projective cover.
By definition of projective objects, the classes of projective maps and surjections are weakly orthogonal and thus define a weak factorization system.
\end{proof}


\medskip

Two maps $u:A\to B$ and $f:X\to Y$ of an arbitrary category $\cC$ are said to be {\it weakly orthogonal} if the pullback hom map $\pbh u f : \Map B Y \to \Map A X \times_{\Map A Y} \Map B Y $ is surjective in $\cS$.
The two classes of projective and surjective maps are weakly orthogonal to each other and form a weak factorization system on $\cS$.
The factorization of $f:X\to Y$ is given by $X\to X+E\to Y$ for $E\to Y$ any surjective map from a set.
This factorization cannot be made functorial.

\medskip
Let $\Delta$ be the category of simplices.
We shall denote the colimit functor (also called realization) $\cS^{\Delta\op}\to \cS$ by $X\index \mapsto |X\index|$ and its right adjoint (constant diagram) by $X\mapsto \underline X$.
The latter functor is fully faithful, and $|-|$ presents $\cS$ as a reflective subcategory of $\cS\splx$.
We denote $W$ the class of maps of $\cS^{\Delta\op}$ sent to invertible maps by the colimit functor.
We call them {\it colimit equivalences}.
For an object $X$ in $\cS$, a {\it resolution} of $X$ is defined as a simplicial diagram $X\index$ equipped with a colimit cocone with apex $X$.
Such a cocone is equivalent to a map $X\index \to \underline X$ which is in $W$.

\medskip
The projective--surjective weak factorization system induces a Reedy weak factorization system on $\cS^{\Delta\op}$.
This weak factorization system is the cofibration--trivial~fibration factorization system of \cite[Theorem 4.4]{MazelGee:model}.
The maps in the right class are the (weak) right orthogonal to the maps $\partial\Delta[n]\to \Delta[n]$ ($n\geq0$).
We shall call them {\it hypersurjective} maps.
A hypersurjective map $U\index\to 1$ is the same thing as a {\it hypercovering} in the topos $\cS$ in the sense of \cite[Definition 6.5.3.2]{Lurie:HTT}.
The following lemma is a crucial property of hypersurjective maps.

\begin{lemma}[{\cite[Proposition 7.2]{MazelGee:model}}]
\label{lem:hyper-colim}
All hypersurjective maps are colimit equivalences.
\end{lemma}
\begin{proof}
This is the proof of $I_{KQ}\mathrm{-inj}\subset \mathbf{W}_{KQ}$ of \cite[Proposition 7.2]{MazelGee:model}.
When the domain is a terminal object, the reader can also look at \cite[Lemma 6.5.3.11]{Lurie:HTT}.
\end{proof}

\noindent The maps in the left class are sometimes called {\it cofibrations}, we shall call them {\it hyperprojective} maps.
A map $X\index\to Y\index$ is hyperprojective if and only if all relative latching maps ${X_n\coprod_{L_nX\index}L_nY\index \to Y_n}$ are projective.
Intuitively, this means that $Y\index$ is build from $X\index$ by adding a {\it set} (rather than an arbitrary space) of non-degenerate simplices in each dimension.
In particular, a map $\emptyset\to X\index$ is hyperprojective if and only if $X\index$ is a simplicial {\it set}.
Thus, the hyperprojective--hypersurjective factorization of $\emptyset\to X\index$ always goes through a simplicial set.
More generally, a map $X\index \to Y\index$ where $X\index$ is a simplicial set is hyperprojective if and only if it is a monomorphism of simplicial sets.
Let $\underline X$ be a constant simplicial object.
Since sets are the projective object of $\cS$, we shall say that a factorization of $\emptyset\to X\index \to \underline X$ is a {\it projective resolution} of $X$.

\medskip
We recall some results on how to construct projective resolutions.
We shall need this to prove \Cref{prop:resolution}.
The following lemma is \cite[Proposition A.2.9.14, Corollary A.2.9.15 and Remark A.2.9.16]{Lurie:HTT}.

\begin{lemma}[Reedy induction]
\label{lem:Reedy-induction}
Let $\cC$ be a category with finite limits and colimits.
The extension of a functor $X:\Delta_{<n}\op \to \cC$ into a functor $X':\Delta_{\leq n}\op \to \cC$ is equivalent to a factorization of the map $L_nX \to M_nX$ (where $L_nX$ and $M_nX$ are the latching and matching objects of $X:\Delta_{<n}\op \to \cC$).
\end{lemma}

\noindent We can apply this to $\cC=\cS$ with the projective--surjective factorization
to get (as a special case of \cite[Corollary 1.4.11]{Lurie:DAGX} or \cite[Corollary 12.4.3.3]{Lurie:SAG}).
\begin{lemma}[Projective resolution]
\label{lem:Reedy-induction2}
Let $X$ be an object in $\cS$.
There exists a simplicial object $\Delta\op \to \cS\slice X$ such that, for every $n$, 
the map $L_nX \to X_n$ is a sum with some set,
and the map $X_n \to  M_nX$ is surjective.
\end{lemma}

\noindent The simplicial object in $\cC\slice X$ of \Cref{lem:Reedy-induction2} provides a map $X\index \to \underline X$ which, by construction, is hypersurjective.

\subsection{Kan groupoids}

This section proves \cref{prop:resolution}, which is going to be our main tool to prove that $\pif\cS$ is locally cartesian closed (\cref{thm:CohLCC}).

\medskip
The following definition is a simplicial analogue of \cite[Definition A.6.5.1]{Lurie:SAG}.
We say that a simplicial space $X\index$ is a {\it Kan groupoid} if it is weakly right orthogonal to all horn inclusions $\Lambda^n_k\to \Delta^n$ ($0\leq k\leq n$, $n\geq 1$), that is if all maps of spaces $X_n \to \Map {\Lambda^n_k}{X\index}$ are surjective in $\cS$.
A simplicial set is a Kan groupoid if and only if it is a Kan complex.
We shall keep the name Kan complex for a Kan groupoid whose values are sets.

The name Kan groupoid is chosen to echo that of Segal groupoid (see below).
In the same way that Segal groupoids can be defined in any category with fiber products, Kan groupoids can be defined in any category with fiber products and a notion of surjection (e.g. a pretopos).
Intuitively, the former provide a notion of ``1-groupoidal equivalence relation'' whereas the later provide a notion of an ``\oo groupoidal equivalence relation''.

\begin{remark}
\label{rem:Kan-fib}
Following \cite[Definition 4.1]{MazelGee:model}, it is convenient to introduce a second weak factorization system on $\cS\splx$, generated by the horn inclusions.
The maps in the right class are called {\it Kan fibrations}.
A simplicial space $X\index$ is a Kan groupoid if and only if $X\index \to 1$ is a Kan fibration.
Hence, if $X\index \to Y\index$ is a Kan fibration and $Y\index$ is a Kan groupoid, then so is $X\index$.
\end{remark}

\begin{lemma}[Kan resolution]
\label{lem:hyper-Kan}
Let $X$ be a space and $\emptyset \to X\index \to \underline X$ be a hyperprojective--hypersurjective factorization of $\emptyset \to X$.
Then the simplicial space $X\index$ is a simplicial set which is a Kan complex.
\end{lemma}
\begin{proof}
The object $X\index$ is a simplicial set by construction of the Reedy factorization system.
We need to prove that it is a Kan groupoid.
Any hypersurjective map is a Kan fibration and by \cref{rem:Kan-fib} it is enough to show that any constant simplicial space $\underline X$ is a Kan groupoid.
Let $|Y\index|$ be the colimit of some simplicial space $Y\index$.
By adjunction, we have $\Map {Y\index} {\underline X} = \Map {|Y\index|} X$.
Applying this to $Y\index=\Lambda^n$ and $Y\index=\Delta^n$ and using that $|\Lambda^n|\to |\Delta^n|$ is an equivalence in $\cS$, we get that $\underline X$ has the lifting condition with respect to all horn inclusions.
\end{proof}

A simplicial space is {\it $n$-coskeletal} if it is the right Kan extension of its restriction to $\Delta_{\leq n} \subset \Delta$.
This is equivalent to the condition that the maps $X_k\to \Map {sk_n \Delta^k} {X_{\leq n}}$ be all equivalence for $k>n$.
We say that a Kan groupoid (in $\cS$) is {\it truncated} if its colimit is $n$-truncated for some $n$.
We say that a Kan complex is has {\it finite values} if is it in $(\Fin\Set)\splx \subset \Set\splx$.

\begin{proposition}
\label{prop:resolution}
A space is unbounded $\pi$\=/finite if and only if it is the geometric realization of a Kan complex with finite values.
A space is $\pi$\=/finite if and only if it is the geometric realization of a coskeletal Kan complex with finite values.
\end{proposition}
\begin{proof}
Let $X\index$ be a Kan complex.
Recall that $\pi_0(|X\index|)$ is a quotient of $X_0$ and 
and $\pi_n(|X\index|,x)$ is a subquotient of $X_n$.
Hence they are all finite if the $X_n$ are.
This proves that the conditions are sufficient.

To see that they are necessary, we use Reedy induction.
Let $Z$ be an unbounded $\pi$\=/finite space, we use \cref{lem:Reedy-induction} in $\cS\slice Z$ to construct a simplicial object $\Delta\op\to \cS\slice Z$.
First, we choose $X_0\to Z$ a surjection from a set $X_0$.
Because $Z$ is unbounded $\pi$\=/finite, $X_0$ can be chosen finite.
At step $1$, we have $L_1(X_{\leq 0}) = X_0$ and $M_0(X_{\leq 0})= X_0\times_Z X_0$.
The space $X_0\times_X X_0$ is a finite sum of path spaces of $Z$.
Since $Z$ is unbounded $\pi$\=/finite, it has a finite number of connected components and we can put $X_1 := X_0 + X'_1$ where $X'_1$ is a finite set (for example $\coprod_{z\in X_0}\pi_1(Z,z)$).
At step $n$, let $\Map {\partial \Delta^n} {X_{<n}}$ be the set of morphisms in $\Set^{\Delta_{<n}\op}$. 
Since all $X_k$ are finite sets, this is a finite set.
Then we have $M_n(X_{<n}) = \Map {\partial \Delta^n} {X_{<n}}\times_{X^{|\partial \Delta^n|}}X$ and $M_n(X_{<n})$ is a sum of $n$-fold path spaces of $Z$.
Since $Z$ is unbounded $\pi$\=/finite, $M_n(X_{<n})$ has a finite number of connected components and we can put $X_n := L_n(X_{<n}) + X'_n$ where $X'_n$ is a finite set (for example $\coprod_{z\in X_0}\pi_n(Z,z)$).
By induction, $L_n(X_{<n})$ is a finite set, hence so is $X_n$.
The resulting simplicial set $X\index$ has finite values.
We get a map $X\index \to \underline Z$ in $\cS\splx$ which is a hypersurjective, hence a colimit cone by \cref{lem:hyper-colim}.
The fact that it is Kan is \cref{lem:hyper-Kan}.
This proves the first assertion.

To prove the second one, we need to recall that for simplicial set $X\index$ with geometric realization $X$, the geometric realization of the morphism $X\index \to cosk_{n+1}X\index$ is the $n$-truncation of $X$.
We leave to reader the proof that $cosk_{n+1}X\index$ is a Kan complex with finite values if $X\index$ is.
By the first part of the proof, the realization of a coskeletal Kan complex with finite values is a $\pi$\-/finite space.
Conversely, if $X$ is a $\pi$\-/finite space which is the realization of a Kan complex $X\index$ with finite values, we can always replace $X\index$ by $cosk_{n+1}X\index$ for some sufficiently large $n$.
This proves the second statement.
\end{proof}

\subsection{Segal groupoids}
\label{sec:groupoids}

This section proves \cref{prop-coh:3}.
We prove in fact a stronger result presenting $\pif\cS$ as the closure of finite sets under quotients of Segal groupoids (\Cref{prop:Segal-gpd}).

\medskip

We say that a simplicial space $X\index$ is a {\it Segal groupoid} if it is right orthogonal to all horn inclusions $\Lambda^n\to \Delta^n$, that is if all maps of spaces $X_n \to \Map {\Lambda^n}{X\index}$ are invertible in $\cS$.

\begin{remark}
\label{rem:Lurie-groupoid}
This definition is equivalent to that of \cite[Definition 6.1.2.7]{Lurie:HTT}.
Using the notations of Proposition 6.1.2.6 in {\it op. cit.}: when $\cC$ has finite limits, the functor $U[\Lambda^k[n]]$ is representable by an object $U^{\Lambda^k[n]}$ in $\cC$.
Then Condition~(3) is equivalent to the canonical map $U_n\to U^{\Lambda^k[n]}$ being invertible, which is our definition of Segal groupoid.
\end{remark}

\begin{definition}[Quotient and effectivity of groupoids]
\label{def:eff-gpd}
Let $X\index$ be a simplicial space and $|X\index|$ its colimit.
We shall also call $|X\index|$ the {\it quotient} of $X\index$ and refer to the canonical map $q:X_0\to |X\index|$ as the {\it quotient map}.
Let $f:X\to Y$ be a map in $\cS$ and $X\index$ be its nerve $N(f)\index$ ($X_n = X\times_Y \dots \times_YX$).
Then, $X\index$ is a Segal groupoid.
Intuitively, $N(f)\index$ is the groupoid encoding the equivalence relation ``to have same image by $f$''.
A Segal groupoid is {\it effective} if the canonical map $X\index \to N(q)$ (where $q$ is the quotient map) is invertible in $\cS\splx$.
In $\cS$, all Segal groupoids are effective, this is part of the Giraud axioms of \oo topoi \cite[Proposition 6.1.3.19]{Lurie:HTT}.
\end{definition}

\begin{remark}
When a simplicial diagram $X\index$ is a Segal groupoid, its quotient is essentially its completion as a Segal space. 
Recall that a complete Segal groupoid is a constant simplicial diagram, then the completion of the Segal groupoid $X\index$ is the constant simplicial diagram with value the quotient $|X\index|$.	
\end{remark}

\begin{proposition}
\label{prop:Segal}
Let $X\index$ be a Segal groupoid in $\pif\cS$, then its quotient $|X\index|$ is in $\pif\cS$.
\end{proposition}
\begin{proof}
The quotient map $X_0\to |X\index|$ being surjective, $\pi_0(|X\index|)$ is a finite set.
Hence we can restrict to the case where $|X\index|$ is connected.
By effectivity of Segal groupoids, we have a cartesian square 
\[
\begin{tikzcd}
X_1 \ar[r]\ar[d] \pbmark & X_0 \ar[d]\\
X_0 \ar[r] & {|X\index|}
\end{tikzcd}
\]
Let $x$ be an element in $X_0$.
The fiber of $X_1\to X_0$ at $x$ is an unbounded $\pi$\=/finite space $Z$ by \cref{prop:342}.
Hence we can apply \cref{prop:342} again to the cartesian square 
\[
\begin{tikzcd}
Z \ar[r]\ar[d] \pbmark & X_0 \ar[d]\\
1 \ar[r] & {|X\index|}
\end{tikzcd}
\]
to deduce that $|X\index|$ is $\pi$\=/finite.
\end{proof}

Lurie proves a similar result for $\Coh\cS$ and Kan groupoids.
We mention it for a comparison.

\begin{proposition}[{\cite[Theorem A.5.5.1]{Lurie:SAG}}]
\label{prop:kan-coh}
Let $X\index$ be a Kan groupoid in $\Coh\cS$, then its quotient $|X\index|$ is in $\Coh\cS$.
\end{proposition}

\medskip
The following result gives meaning to the category $\pif\cS$ and $\Coh\cS$ inside $\cS$.
We shall use it in \cref{thm:initial} to prove $\pif\cS$ is the initial \oo pretopos.

\begin{proposition}[Exact completions]
\label{prop:Segal-gpd}
\begin{enumerate}
\item \label{prop:Segal-gpd:2} The category $\pif\cS$ is the smallest subcategory of $\cS$ containing $\Fin\Set$ and closed under quotients of Segal groupoids.    
\item \label{prop:Segal-gpd:1} The category $\Coh\cS$ is the smallest subcategory of $\cS$ containing $\Fin\Set$ and closed under quotients of Kan groupoids.
\end{enumerate}
\end{proposition}
\begin{proof}
\noindent \eqref{prop:Segal-gpd:2}
Let $\cC\subset \cS$ be the smallest full subcategory containing $\Fin\Set$ and closed under quotients of Segal groupoids.
Since $\Fin\Set\subset \pif\cS$, \Cref{prop:Segal} proves that $\cC\subset \pif\cS$.
Conversely, we proceed by induction on the truncation level.
Let $\pif\cS\truncated n\subset \pif\cS$ be the full subcategory spanned by $n$-truncated objects.
We have $\pif\cS\truncated 0 = \Fin\Set$.
Let us prove that any object $X$ of $\pif\cS\truncated {n+1}$ can be obtained as the quotient of a Segal groupoid in $\Coh\cS\truncated n$.
Let $f:X_0\to X$ be a surjective map where $X_0$ is a set, then $X_0$ is in $\pif\cS\truncated n$.
We consider the nerve $X\index$ of $f$.
It is a Segal groupoid whose quotient is $X$.
The result will be proved if we show that $X\index$ is a simplicial object in $\pif\cS\truncated n$.
The space $X_1 = X_0 \times_X X_0$ is a sum of loop spaces of $X$, hence $\pi$\=/finite and $n$-truncated.
More generally, we have $X_n = X_1\times_{X_0}\dots \times_{X_0}X_1$, and this shows $X_n$ is also in $\pif\cS\truncated n$.

\smallskip
\noindent \eqref{prop:Segal-gpd:1}
By \cref{prop:resolution}, $\Coh\cS$ is included in the smallest full subcategory containing $\Fin\Set$ and closed under quotients of Kan groupoids.
The converse is given by \cref{prop:kan-coh}.
\end{proof}

\begin{corollary}[Descent properties]
\label{prop:Segal-descent}
\begin{enumerate}
\item \label{prop:Segal-descent:1} Quotients of Kan groupoids have descent in $\Coh\cS$.
\item \label{prop:Segal-descent:2} Quotients of truncated Kan groupoids have descent in $\pif\cS$.
\item \label{prop:Segal-descent:3} Segal groupoids have descent in $\pif\cS$ (\cref{def:descent}).
\item \label{prop:Segal-descent:4} Segal groupoids are universal and effective in $\pif\cS$ (\cref{def:eff-gpd}).
\end{enumerate}
\end{corollary}
\begin{proof}
\noindent The properties
\eqref{prop:Segal-descent:1},
\eqref{prop:Segal-descent:2}, and
\eqref{prop:Segal-descent:3} 
are consequences of \cref{lem:descent}.
We are left to prove \eqref{prop:Segal-descent:4}.
The universality of Segal groupoids is a consequence of \eqref{prop:Segal-descent:3},
and the effectivity is a consequence of $\pif\cS\subset \cS$ preserving finite limits (\cref{prop:coh-lex}).
\end{proof}

\begin{remark}
Putting together \Cref{prop:coh-lex,prop:coh-extensive,rem:coh,prop:Segal-descent}, we get that $\Coh\cS$ and $\pif\cS$ are \oo pretopoi in the sense of \cite[Definition A.6.1.1]{Lurie:SAG}.
This can also be deduced by applying \cite[Corollary A.6.1.7]{Lurie:SAG} to $\cS$.
We shall see in \cref{thm:initial} that $\pif\cS$ is the initial pretopos.
Since $\Fin\Set$ is the initial 1-pretopos, $\pif\cS$ is then the \oo pretopos envelope of $\Fin\Set$, and can be thought as its higher exact completion.
\end{remark}

\subsection{Local cartesian closure}

This section proves \Cref{prop-coh:5} (\cref{thm:CohLCC}).
We prove first that $\pif\cS$ is cartesian closed and deduce the statement for the slice categories by a descent argument.

\medskip
We start with some technical lemmas.
We say that a functor $f:C\to D$ between $n$-categories is {\it $n$-final} if for any cocomplete $n$-category $\cC$, the colimit of any diagram $X:D\to C$ coincide with the colimit of $X\circ f:C\to \cC$.
Dually, a functor is {\it $n$-initial} if $f\op:C\op\to D\op$ is $n$-final.
For $C$ a small $n$-category, its free cocompletion (as an $n$-category) is 
$\relP n C := \fun {C\op} {\cS \truncated {n-1}}$.

\begin{lemma}
\label{lem:final-general}
The following condition are equivalent:
\begin{enumerate}
\item \label{lem:final-general:1} the functor $f:C\to D$ is $n$-final;
\item \label{lem:final-general:2} the functor $\relP n f : \relP n C \to \relP n D$ preserves the terminal object;
\item \label{lem:final-general:3} for any $d$ in $D$, the realization of the category $C\coslice d := C\times_DD\coslice d$ is an $(n-1)$-connected space.
\end{enumerate}
\end{lemma}

\begin{proof}
\eqref{lem:final-general:1} $\Rightarrow$ \eqref{lem:final-general:2}
The colimit of the Yoneda embedding $C\to \relP n C$ is the terminal object.
If $\cC$ is a cocomplete $n$-category the colimit of a diagram $C\to \cC$ is the image of the terminal object by the left Kan extension $\relP n C \to \cC$.

\noindent \eqref{lem:final-general:2} $\Rightarrow$ \eqref{lem:final-general:1}
follows from the fact that, given a diagram $D\to \cC$, we have a commutative diagram
\[
\begin{tikzcd}
C \ar[d,"f"']\ar[r] & \relP n C\ar[d,"\relP n f"] \ar[rd,dashed, bend left]\\  
D \ar[r] & \relP n D \ar[r, dashed]& \cC
\end{tikzcd}
\]
(where the dashed arrows are left Kan extensions along the Yoneda embeddings).

\smallskip
\noindent\eqref{lem:final-general:2} $\Leftrightarrow$ \eqref{lem:final-general:3}
We have $\relP n f(1) = 1$ if and only if, for any $d$ in $D$, the space $\Map d {\relP n f(1)}$ is contractible.
But, in $\cS\truncated {n-1}$, we have
\[
\Map d {\relP n f(1)} =  \Map {\widehat d} {\colim_C \widehat{f(c)}} = \colim_C \Map d {f(c)} = \colim_{C\coslice d}1 = |C\coslice d|\truncated {n-1}
\]
where $|C\coslice d|\truncated {n-1}$ is the $(n-1)$-truncation of the realization of $C\coslice d$.
This space is contractible if and only if the realization of $C\coslice d$ is $(n-1)$-connected.
This proves \eqref{lem:final-general:2} $\Leftrightarrow$ \eqref{lem:final-general:3}.
\end{proof}

\begin{lemma}
\label{lem:final}
The inclusion $\Delta_{\leq n}\to \Delta$ is initial for diagrams in $n$-categories.
\end{lemma}
\begin{proof}
This can be seen as an application of \cite[Lemma 6.5.3.10]{Lurie:HTT}.
We will give a more direct informal argument.
If $k\leq n$, $(\Delta_{\leq n})\slice {[k]}$ has a terminal object and is weakly contractible.
If $k>n$, the  realization of $(\Delta_{\leq n})\slice {[k]}$ is $sk_n(\Delta[k])$ which is a bouquet of $n$-spheres, hence $(n-1)$-connected.
Then the result follows from the dual of \Cref{lem:final-general}~\eqref{lem:final-general:3}.
\end{proof}

Recall that an object $X$ in category $\cC$ is called $n$-truncated if the functor $\Map - X:\cC\op\to \cS$ takes values in $n$-truncated spaces. 

\begin{lemma}
\label{lem:final2}
The inclusion $\Delta_{\leq n}\to \Delta$ is initial for diagrams of $(n-1)$-truncated objects.
\end{lemma}
\begin{proof}
Let $\cC$ be a category and $\cC\truncated {n-1}\subset \cC$ be the full subcategory of $n$-truncated objects.
Let $D:I\to \cC\truncated {n-1}$ be a diagram having a limit in $\cC$. 
Let us see that its limit is in $\cC\truncated n$.
The result is true in $\cS$ because the subcategory $\cS\truncated {n-1}\subset \cS$ of $n$-truncated spaces is reflective, hence closed under arbitrary limits.
For a general $\cC$, the limit of $D$ is the object representing the functor 
\begin{align*}
\cC\op &\tto \cS   \\
X &\mto \lim_i\Map X {D_i}
\end{align*}
If all the $D_i$ are $(n-1)$-truncated, this functor takes values in $(n-1)$-truncated spaces, so any representative will be an $(n-1)$-truncated object.
This reduces the problem to prove the initiality of $\Delta_{\leq n}\to \Delta$ to diagrams in the $n$-category $\cC \truncated {n-1}$, but then it follows from \Cref{lem:final}.
\end{proof}

\begin{proposition}
\label{prop:CC}
The category $\pif\cS$ is cartesian closed, and the embedding $\pif\cS\subset \cS$ preserves the exponentials.
\end{proposition}
\begin{proof}
We have seen that $\pif\cS\subset \cS$ is closed under finite products (\cref{prop:coh-lex}).
We are going to show that for any two spaces $X$ and $Y$ in $\pif\cS$, the space $Y^X$ is in $\pif\cS$.
When $X$ is a finite set, this is true because $Y^X$ is a finite product of $Y$.
For a general $X$, we use \Cref{prop:resolution} to present $X$ as a colimit of a simplicial finite set, and get
\[
Y^X = \lim_m Y^{X_m}.
\]
This limit is a priori infinite and $\pif\cS$ is only closed under finite limits.
By assumption $Y$ is $k$-truncated for some $k$, then so are all the $Y^{X_m}$.
Thus, we can use \cref{lem:final2} and replace the limit by an equivalent one which is finite. 
\end{proof}

\begin{lemma}
\label{lem:limCCC}
The limit of a diagram of cartesian closed categories and cartesian closed functors is cartesian closed.
\end{lemma}
\begin{proof}
The proof is straightforward for products, so we need only to give an argument for fiber products.
Let $\cC_1\xto p \cC_0 \xot q \cC_2$ be a diagram of cartesian closed categories and cartesian closed functors.
The objects in the limit $\cC:=\lim \cC_i$ are families $X = (X_1,X_0,X_2,x_1:p(X_1)\simeq X_0,x_2:q(X_2)\simeq X_0)$.
We claim that the internal hom between two such families $X$ and $Y$ are computed termwise as
\[
X^Y := \big(
X_1^{Y_1},\ 
X_0^{Y_0},\ 
X_2^{Y_2},\ 
x_1^{y_1^{-1}}\!:\,p(X_1)^{p(Y_1)}\simeq X_0^{Y_0},\ 
x_2^{y_2^{-1}}\!:\,q(X_2)^{q(Y_2)}\simeq X_0^{Y_0}
\big)\,.
\]
To verify it, we will build a natural equivalence $\Map{Z\times Y}X \simeq \Map{Z}{X^Y}$ for $Z$ a third object of $\cC$.
The mapping space $\Map{Z\times Y}X$ is given by
\begin{equation}
\label{eq:limCCC:prod}
\Map{Z_1\times Y_1}{X_1}
\underset{\Map{p(Z_1)\times p(Y_1)}{X_0}}{\times}
\Map{Z_0\times Y_0}{X_0}
\underset{\Map{q(Z_2)\times q(Y_2)}{X_0}}{\times}
\Map{Z_2\times Y_2}{X_2}\,,
\end{equation}
and $\Map{Z}{X^Y}$ by
\begin{equation}
\label{eq:limCCC:exp}
\Map{Z_1}{X_1^{Y_1}}
\underset{\Map{p(Z_1)}{X_0^{Y_0}}}{\times}
\Map{Z_0}{X_0^{Y_0}}
\underset{\Map{Z_2}{X_2^{Y_0}}}{\times}
\Map{q(Z_2)}{X_2^{q(Y_2)}}\,.
\end{equation}
In each $\cC_i$, we have a natural ``transposition'' equivalence $\lambda^2:\Map{Z_i\times Y_i}{X_i}\simeq\Map{Z_i}{X_i^{Y_i}}$.
Using these equivalences and the fact that $p$ preserves products and internal homs, we leave the reader check that we can get a commutative diagram
\[
\begin{tikzcd}
\Map{Z_1\times Y_1}{X_1}
	\ar[rr,"x_1\circ p(-)"] \ar[dd,"\simeq"',"\lambda^2"]
&& \Map{p(Z_1)\times p(Y_1)}{X_0}
	\ar[from=rr,,"-\circ (z_1\times y_1)"'] \ar[dd,"\simeq"',"X_0^{y_1^{-1}}\circ \lambda^2"]
&& \Map{Z_0\times Y_0}{X_0}
	\ar[dd,"\simeq"',"\lambda^2"]
\\
\\
\Map{Z_1}{X_1^{Y_1}}
	\ar[rr,"\left(x_1^{y_1^{-1}}\right)\circ p(-)"]
&& \Map{p(Z_1)}{X_0^{Y_0}}
	\ar[from=rr,"-\circ z_1"'] 
&& \Map{Z_0}{X_0^{Y_0}}
\end{tikzcd}
\]
whose vertical maps are equivalences.
By taking the horizontal limits, we get a natural equivalence between the left fiber products of \eqref{eq:limCCC:prod} and \eqref{eq:limCCC:exp}.
Proceeding the same way for the right side, we get the expected natural equivalence $\Map{Z\times Y}X \simeq \Map{Z}{X^Y}$.
\end{proof}

\begin{remark}
A more conceptual proof of \cref{lem:limCCC} is that the category of cartesian closed categories is monadic over that of categories, hence the forgetful functor creates limits.
\end{remark}

\begin{theorem}
\label{thm:CohLCC}
The category $\pif\cS$ is locally cartesian closed.
Moreover, for each $X$ in $\pif\cS$ the embedding $(\pif\cS)\slice X\subset\cS\slice X$ preserves the internal hom.
\end{theorem}

\begin{proof}
We need to prove that for any $X$ in $\pif\cS$, the category $(\pif\cS)\slice X$ is cartesian closed.
We use \cref{prop:resolution} to get a Kan complex with finite values $X\index$ with colimit $X$.
\Cref{prop:Segal-descent} gives that $(\pif\cS)\slice X = \lim_{\Delta}\ (\pif\cS)\slice {X_n} = \lim_{\Delta}\prod_{X_n} \pif\cS$ and the result follow from \cref{lem:limCCC}.

For the second statement, the case $X=1$ is \cref{prop:CC}.
For, the general case, we choose a surjective family of points $x_i:1\to X$.
Then the families of pullbacks $x_i^*:(\pif\cS)\slice X\to \pif\cS$ and $x_i^*:\cS\slice X\to \cS$ are conservative and preserve internal homs.
The inclusion $\iota:(\pif\cS)\slice X\subset\cS\slice X$ preserves internal hom 
if and only if, for every $i$, the induced functors $(\pif\cS)\to \cS$ do.
Since this last functor is equivalent to the case $X=1$, this finishes the proof.
\end{proof}

\begin{remark}
This result is not true for the category $\Coh\cS$.
Let $X = \prod_n K(\ZZ_2,n)$.
Any sequence of group morphisms $\phi_n: \ZZ_2\to \ZZ_2$ defines an endomorphisms $\phi:=\prod_n K(\phi_n,n)$ of $X$.
Acting differently on the $\pi_n$, these $\phi$ are non-homotopic in $X^X$.
Any group morphism $\ZZ_2\to \ZZ_2$ is either the identity or constant.
Hence, the set of such sequences is $2^\NN$. 
This proves $\pi_0(X^X)$ is not finite (it's not even countable).
\end{remark}

\section{The universe of \texorpdfstring{$\pi$}{pi}-finite spaces}

This section proves \cref{prop-coh:6,prop-coh:7,prop-coh:8}.
We do so by constructing first a universe for $\pi$\=/finite spaces in $\cS$ (\cref{prop:pi-finite-maps,thm:universe}).
We also prove the negative \cref{fact:2} in \cref{thm:no-universe}.
We will rely heavily on the material of \cref{app:sec:univ-maps,app:sec:univ-refl,app:sec:univ-topos,app:sec:univ-structure} and the reader is advised to read it first.

\subsection{The universe of \texorpdfstring{$\pif\cS$}{pi-finite spaces} in \texorpdfstring{$\cS$}{spaces}}

The topos $\cS$ is localic thus a bounded topos \cite[Example A.7.1.3]{Lurie:SAG}.
It is also coherent and locally coherent \cite[Example A.2.1.7]{Lurie:SAG}.
The subcategory $\pif\cS\subset\cS$ is the associated pretopos of bounded coherent objects and we have therefore $\cS=\Sh{\pif\cS}$ by \cite[Theorem A.7.5.3]{Lurie:SAG} (where the sheaves are taken for the coherent topology, see \cref{def:coh-topo}).
The {\it universe} of $\pif\cS$ is the functor $\pif\sfS:\pif\cS\op\to \cS$ sending $X$ to $(\pif\cS)\slice X\int$, the internal groupoid of the slice category $(\pif\cS)\slice X$ (see \cref{def:universe}).
By \cref{prop:descent-pretopos}, it is a sheaf for the coherent topology, and thus defines an object of $\cS$, that we denote $\pif U$.
We denote $\pif u:\pif U'\to \pif U$ the associated universal family of $\pif \cS$ (\cref{app:def:univ-fam}).
By \cref{app:prop:univ-fam=univ}, it is a univalent map in $\cS$.
By \cref{app:sec:univ-fam-S}, $\pif u$ can be described as a sum
\[
\pif u : U'_\pi\stto U_\pi
\quad:=\quad 
\coprod_{F\in\Sigma} F/Aut(F)
\stto \coprod_{F\in\Sigma} BAut(F)
\]
where $F$ runs over a set $\Sigma$ of representative for the isomorphism classes of the fibers of $\pif u$, which are all the $\pi$\=/finite spaces by definition of $\pif\sfS$.
The following lemma shows that $\Sigma$ is a countable set.

\begin{lemma}
\label{lem:coh-countable}
The set of isomorphism classes of objects of $\pif\cS$ is countable.
\end{lemma}
\begin{proof}
It is enough to show the result for $n$-truncated $\pi$\=/finite space.
\Cref{prop:resolution} proves that every such $\pi$\=/finite space $X$ can be described as the colimit of a simplicial diagram $X\index:\Delta\op\to \Fin\Set$ which is $n$-coskeletal.
This shows that every $n$-truncated $\pi$\=/finite space can be presented by means of a diagram $\Delta_{\leq n+1}\op\to \Fin\Set$.
The result follows since the set of isomorphism classes of such diagrams is countable.
\end{proof}

\begin{remark}
\label{lem:countable}
Notice that the connected component of $\pif U$ are $\pi$\=/finite spaces, and therefore realization of a simplicial set with finite value (\cref{prop:resolution}).
Then \cref{lem:coh-countable} implies that $\pif U$ is the realization of a countable simplicial set.
\end{remark}

We consider the map $\pif f:=\coprod_{\Sigma} F \to \coprod_{\Sigma} 1$.
There is an obvious cartesian surjection $\pif f\to \pif u$ and therefore $\pif u$ is the univalent reflection of $\pif f$ by \cref{app:lem:equiv-univ-reflection}.
In particular, we have an equality of the local classes of maps $\{\pif f\}\loc=\{\pif u\}\loc$ by \cref{lem:fam-equiv}.
We shall say that a map in $\Sh\cE$ is a {\it $\pi$\=/finite map} if all its fibers are $\pi$\=/finite spaces.

\begin{proposition}
\label{prop:pi-finite-maps}
The local class of maps $\{\pif u\}\loc$ is the class of $\pi$\=/finite maps.
\end{proposition}
\begin{proof}
This follows from the discussion in \cref{app:sec:univ-fam-S}, but we reproduce the argument.
Since $\pif u$ is a univalent map, it is terminal in $\{\pif u\}\loc$ and a map is in $\{\pif u\}\loc$ if and only if it is a pullback of $\pif u$.
In particular, the fibers of such a map are $\pi$\=/finite spaces.
Conversely, if $f:X\to Y$ is a map whose fibers are $\pi$\=/finite spaces, we fix a surjection $b:Y_0\surj Y$ from a set.
The pullback $f'$ of $f$ along $b$ is a sum of maps $F\to 1$ where $F$ is a $\pi$\=/finite space.
In particular, $f'$ is a pullback of the map $\pif f$.
This provides a span $\pif f \ot f' \surj f$ which shows that $f$ belong to $\{\pif f\}\loc=\{\pif u\}\loc$.
\end{proof}

\begin{theorem}
\label{thm:closure-pifU}\label{thm:universe}
The univalent map $\pif u:\pif U'\to \pif U$ is closed under diagonals, dependent sums, dependent products, finite sums, finite products, and quotients of Segal groupoids (see \cref{app:sec:univ-structure} for definitions).
\end{theorem}
\begin{proof}
Diagonals. The diagonal of a $\pi$\=/finite map $f:X\to Y$ is a $\pi$\=/finite map since its fibers are path spaces of $\pi$\=/finite spaces.
Then the statement follows from \cref{prop:closure-diagonal}.

\smallskip
\noindent
Dependent sums. 
By \cref{prop:pi-finite-maps}, the local class $\{\pif u\}\loc$ is that of $\pi$\=/finite maps.
This class is closed under composition by \cref{prop:342} (applied separately on each connected component).
Then the statement follows from \cref{lem:sigma=compo}.

\smallskip
\noindent
Dependent products. By \cref{def:dep-sigma-pi}~\eqref{def:dep-sigma-pi:2} we must see that, for any two maps $f:X\to Y$ and $g:X'\to X$ in $\cM$, the map $f_*(g)$ is in $\cM$.
This is a consequence of \cref{thm:CohLCC}.

\smallskip
\noindent
Finite sums.
Given a finite family $X_i\to Y$ of $\pi$\=/finite maps, the fibers of $\coprod_i X_i\to Y$ are finite sums of $\pi$\=/finite spaces, thus $\pi$\=/finite spaces by \cref{lem:coh-extensive}.
This proves \cref{def:fin-sigma-pi}~\eqref{def:fin-sigma-pi:local}.
The proof for finite products is similar using \cref{prop:coh-lex}.

\smallskip
\noindent
Quotients.
Let $X\index \to Y$ be a Segal groupoid in $\Sh\cE\slice Y$ such that all maps $X_n\to Y$ are $\pi$\=/finite.
By universality of colimits in $\Sh\cE$, the colimit of $X\index$ can be computed fiberwise over $Y$.
By \cref{prop:Segal} applied fiberwise, the quotient map $|X\index|\to Y$ is a $\pi$\=/finite map.
This proves \cref{app:def:Segal}~\eqref{app:def:Segal:local}.
\end{proof}

\begin{remark}
In terms of dependent type theory, \cref{thm:closure-pifU} says that the map $\pif u:\pif U'\to \pif U$ provides a model for a type theory with identity types, $\Sigma$-types, $\Pi$-types,  finite sums, finite products, and quotients of equivalence relations.
\end{remark}

Recall the universal family of finite sets from \cref{app:ex:fin-set}
$u_\set:\coprod_n \underline n/\mathfrak S_n \to \coprod_n B\mathfrak S_n$ (where $\mathfrak S_n$ is the group of permutations of $n$ elements).
Since finite sets are $\pi$\=/spaces, this is a subfamily of $\pif u$.

\begin{corollary}
\label{cor:initial-LCCpretopos}
The family $\pif u$ is the smallest family containing $u_\set$ and closed under quotients of Segal groupoids.
\end{corollary}
\begin{proof}
The family $\pif u$ is closed under quotients of Segal groupoids by \cref{thm:closure-pifU}.
Then the statement is essentially \cref{prop:Segal-gpd}~\eqref{prop:Segal-gpd:2} (applied fiberwise).
We leave the details to the reader.
\end{proof}

\subsection{Univalent families in \texorpdfstring{$\pif\cS$}{pi-finite spaces}}
\label{sec:univ-fam-pifS}

The notion of a pretopos with enough univalent families is defined in \cref{app:sec:univ-refl}.

\begin{theorem}[Enough univalent families]
\label{thm:univalent}
The category $\pif\cS$ has enough univalent families.
Moreover, every univalent map is of the form
\[
u_S:U'_S\stto U_S
\quad:=\quad 
\coprod_{F\in S} F/Aut(F)
\stto \coprod_{F\in S} BAut(F)
\]
where $S$ is a finite set of two by two non-isomorphic $\pi$\=/finite spaces.
\end{theorem}
\begin{proof}
The first statement is an application of \cref{app:main:cor}, since $\pif\cS$ is a locally cartesian closed pretopos by \cref{thm:CohLCC}.
The second statement can be proved in several ways.
For example by a reasoning similar to the description of univalent maps in $\cS$ of \cref{app:sec:univ-fam-S} (using the fact that every $\pi$\=/finite space can be covered by a finite set).
Or by characterizing these maps as the subfamilies of $\pif u$ in $\cS$ which are maps in $\pif\cS$.
\end{proof}

\begin{remark}
The illustration of \cref{app:prop:univ=sum-univ} is particularly easy here:
it is clear that $\pif u$ is the filtered union of all $u_S$.
\end{remark}

\begin{lemma}
\label{lem:Omega-univ}
The inclusion $\top:\{1\}\to \{\emptyset,1\}$ is univalent in $\cS$.
\end{lemma}
\begin{proof}
For $S:= \{\emptyset, 1\}$, we have $Aut(\emptyset) = Aut(1) = 1$ and 
$U_S = BAut(\emptyset) + BAut(1) = 1+1$ and $U'_S = \emptyset/Aut(\emptyset) +1/Aut(\emptyset) = \emptyset+1=1$. 
Thus, the map $\top$ is isomorphic to the univalent map $U'_S\to U_S$.
\end{proof}

The notion of a family closed under diagonals is defined in \cref{app:sec:identity-types}.

\begin{proposition}[Univalent families with diagonals]
\label{thm:univalent:id}
For any map $f:X\to Y$ in $\pif\cS$, there exists a univalent map $u$ in $\pif\cS$, which is closed under diagonals and such that every diagonal of $f$ is a base change of $u$.
\end{proposition}
\begin{proof}
The map $f$ being truncated, there exists an $n$ such that for all $k\geq n$ the iterated diagonal $\Delta^kf$ are invertible.
In particular, for such a $k$, $\Delta^kf$ will always be a pullback of $\Delta^nf$.
We put $g:=\coprod_{0\leq k\leq n} \Delta^kf$.
Using the effectivity of sums, we get that the map $\Delta g$ is a pullback of $g$.
So $g$ is a family closed under diagonals.
By \cref{lem:univ-diagonal}, the univalent reflection $u_g$ of $g$ is also closed under diagonals.
The composition of cartesian morphisms $f\to g\to u_g$ shows that $f$ is a base change of $u_g$.
Then, the result about the diagonals of $f$ follows from \cref{prop:closure-diagonal}.
\end{proof}

\begin{proposition}
\label{prop:CohOmega}
The set $2:=\{\emptyset,1\}$ is a (Boolean) subobject classifier in $\pif\cS$,
and the univalent family $\top:1\to 2$ of \cref{lem:Omega-univ} is the universal subobject.
Moreover, the family $\top$ is closed under dependent sums, dependent products and diagonals.
\end{proposition}
\begin{proof}
The map $\top:1\to 2$ is a universal subobject in $\Fin\Set$.
A map $X\to Y$ in $\cS$ is a monomorphism if and only if the map $\pi_0(X)\to \pi_0(Y)$ is injective.
If $Sub(X)$ is the set of subobjects of a space $X$, we have natural bijections
\[
Sub(X) = 2^{\pi_0(X)} = 2^X.
\]
This proves that 2 is a also subobject classifier in $\cS$.
It is moreover Boolean since all posets $Sub(X)$ are Boolean algebras.
Then, the result in $\pif\cS$ follows from the fact that any subobject of a $\pi$\=/finite space $X$ is $\pi$\=/finite (\cref{lem:coh-subobject}).

The map $\top:1\to 2$ is univalent by \cref{lem:Omega-univ}.
We just saw that the associated local class is that of monomorphisms in $\pif\cS$.
Since monomorphisms are closed under composition, $\top$ has dependent sums by \cref{lem:sigma=compo}.
It has dependent products because the functor $f_*$ being right adjoint, it always preserves monomorphisms (\cref{def:dep-sigma-pi}~\eqref{def:dep-sigma-pi:2}).
The closure under diagonals follows from the fact that the diagonal of a monomorphism is an isomorphism thus a monomorphism.
\end{proof}

\begin{remark}
The map $\top:1\to 2$ is not closed under finite sums and product in the sense of \cref{def:dep-sigma-pi}.
Nonetheless it has finite sums and products in the sense that subobjects in $\pif\cS$ are closed under finite unions and intersections.
\end{remark}

\begin{lemma}
\label{lem:sigma-pi}
Let $u:U'\to U$ be a univalent family in $\cS$ and let $X$ be one of the fibers of $u$.
\begin{enumerate}
\item\label{lem:sigma-pi:1} If the map $u$ is closed under dependent sums, then all finite powers $X^n$ are fibers of $u$.
\item\label{lem:sigma-pi:2} If the map $u$ is closed under dependent products, then all iterated exponential $X^X$, $X^{(X^X)}$, ... are fibers of $u$.
\end{enumerate}
\end{lemma}
\begin{proof}
\eqref{lem:sigma-pi:1}
The map $q:X\to 1$ is a base change of $u$ by assumption.
The map $p_1:X\times X\to X$ is a base change of $q$ and thus of $u$.
The composition $qp_1 = q_!(p_1)$ is then a base change of $u$ since $u$ is closed under dependent sums.
The higher powers are obtained from there by an induction left to the reader.

\smallskip
\noindent \eqref{lem:sigma-pi:2}
We proceed as in \eqref{lem:sigma-pi:1}, the map $q_*(p_1)=X^X\to 1$ is then a base change of $u$ since $u$ is closed under dependent products.
The higher iterated exponential are obtained by induction.
\end{proof}

\begin{theorem}
\label{thm:no-universe}
The map $\top:1\to 2$ is the largest univalent map in $\pif\cS$ with dependent sums and the largest univalent map with dependent products.
\end{theorem}
\begin{proof}
The map $\top:1\to 2$ is univalent by \cref{lem:Omega-univ} and has dependent sums and product by \cref{prop:CohOmega}.
Let us see that it is maximal with these properties.
Let $U'_S\to U_S$ be a univalent map in $\pif\cS$ such that one of the components of $U_S$ is $BAut(X)$ for a $\pi$\=/finite space $X$ which is not $\emptyset$ or 1.
Let us see now that $U'_S\to U_S$ cannot be closed under dependent sums.
If this was the case, by \cref{lem:sigma-pi}~\eqref{lem:sigma-pi:1}, all finite powers $X^n$, would be classified by $U_S$.
When $X$ is not subterminal, it has a non-trivial $\pi_n(X,x)$ for some $n$ and some point $x$.
Therefore all finite powers $X^n$ have non-isomorphic $\pi_n$.
This implies that $U_S$ must have a countable number of connected components, which contradicts the fact that it belong to $\pif\cS$.
This proves that any object $X$ classified by $U_S$ must be subterminal.
The argument is similar if we consider a univalent family with dependent products (using \cref{lem:sigma-pi}~\eqref{lem:sigma-pi:2}).
\end{proof}

\section{Initiality properties}

This section proves the initiality results of \cref{prop-coh:9,prop-coh:10bis,prop-coh:10}.

\subsection{Initial pretopos}

Recall from \cite[Definition A.6.1.1]{Lurie:SAG}, that an \oo pretopos (we shall say simply a {\it pretopos}) is a category $\cE$
with finite limits,
with extensive finite sums,
and with universal and effective quotients of Segal groupoids.
A morphism of pretopoi is a functor preserving finite limits, finite sums and quotients of Segal groupoids (or equivalently surjective maps).

\begin{theorem}
\label{thm:initial}
The category $\pif\cS$ is the initial pretopos.
\end{theorem}

\begin{proof}
Let $\cE$ be a pretopos.
We consider the following inclusions of categories
\[
\fun{\pif\cS}\cE^\mathrm{lex}_{\mathrm{Segal},\sqcup}
\quad\subseteq\quad
\fun{\pif\cS}\cE^1_{\mathrm{Segal},\sqcup}
\quad\subseteq\quad
\fun{1}\cE^1
\quad=\quad 1\,.
\]
The category $\fun{\pif\cS}\cE^\mathrm{lex}_{\mathrm{Segal},\sqcup}$ is that of morphisms of pretopoi (preserving finite sums, quotient of Segal groupoids and  finite limits).
It is a full subcategory of $\fun{\pif\cS}\cE^1_{\mathrm{Segal},\sqcup}$ (functors preserving finite sums, quotient of Segal groupoids and the terminal object).
By \cref{prop:Segal-gpd}, the terminal object of $\pif\cS$ is dense in $\pif\cS$ and every object can be accessed by finite sums and quotient of Segal groupoids.
This provide a fully faithful functor $\fun{\pif\cS}\cE^1_{\mathrm{Segal},\sqcup}\subseteq \fun{1}\cE^1$ into the category of terminal objects of $\cE$.
Since the latter category is contractible, this shows that $\fun{\pif\cS}\cE^\mathrm{lex}_{\mathrm{Segal},\sqcup}$ is either empty of contractible.

The proof will be finished if we produce a pretopos morphism $\pif\cS\to\cE$.
Let $\Sh\cE$ be the topos of sheaves on $\cE$ for the effective epimorphism topology in the sense of \cite[Definition A.6.2.4]{Lurie:SAG}.
The topos $\cS$ is initial in the category of topoi and cocontinuous and left-exact functors (the opposite category of topoi and geometric morphisms).
Let $i:\cS\to \Sh\cE$ be the unique such functor.
By definition of the topology, the inclusion $\cE\subseteq\Sh\cE$ preserves finite sums, quotient of Segal groupoids and finite limits.
The functor $i$ send the terminal object of $\cS$ to that of $\Sh\cE$, which is also that of $\cE$.
Since it preserve fintie sums, it sends $\set\subset\cS$ to $\cE$.
Then, using \cref{prop:Segal-gpd}, the image of $\pif\cS\to\Sh\cE$ must then be in $\cE$.
This proves that $i$ restricts to a pretopos morphism $\pif\cS\to \cE$.
\end{proof}

\subsection{Initial \texorpdfstring{$\Pi$}{Pi}-pretopos}
\label{sec:pi-pretopos}

Let $\cE$ be a pretopos and $i:\pif\cS\to \cE$ the morphism of \cref{thm:initial}.
For $X$ a space, we denote by $\cE^X$ the category of $X$-diagrams in $\cE$. 
We show it is equivalent to $\cE\slice {iX}$.

\begin{lemma}
\label{lem:bcoh-diagrams}
For $X$ a $\pi$\=/finite space, there exists a canonical equivalence $\cE\slice {iX}=\cE^X$.
Moreover, given a morphism of pretopoi $f:\cE\to \cF$, the pretopos morphism $f\slice X:\cE\slice {iX}\to \cF\slice {iX}$ corresponds under this equivalence to the pretopos morphism $f^X:\cE^X\to \cF^X$.
\end{lemma}
\begin{proof}
We prove it by descent.
When $X$ is a finite set this is true 
by extensionality of sums in $\cE$ and because $i:\pif\cS\to \cE$ preserves finite sums.
For a general $X$, we use a resolution $X\index$ by a truncated Kan complex (\cref{prop:resolution}).
By \cref{thm:initial} we have $iX = \colim i(X_n)$ in $\cE$ (we shall simply write $X_n$ for $i(X_n)$ henceforth).
By the descent property of \cref{prop:Segal-descent} and extensivity, we get $\cE\slice {iX} = \lim_n \cE\slice {X_n}= \lim_n \cE^{X_n}$.
Recall that the embedding $\cS\subset \Cat$ of groupoids in categories preserves all limits and colimits (since it has both a left and a right adjoint).
This gives $\lim_n \cE^{X_n} = \cE^{\colim X_n} = \cE^X$.
Altogether, this provides the equivalence of the first statement.

For the second statement, we denote by $i:\pif\cS\to \cE$ and $j:\pif\cS\to \cF$ the canonical morphisms of \cref{thm:initial}. 
We have $fi=j$.
We need to show that there exists a natural commutative square
\[
\begin{tikzcd}
\cE\slice {iX} \ar[d,"f\slice X"']\ar[r,"\simeq"] & \cE^X \ar[d,"f^X"]\\
\cF\slice {jX} \ar[r,"\simeq"] & \cF^X\,.
\end{tikzcd}
\]
The top-then-down composition sends a map $Y\to X$ to the family $x:1\to X\mapsto f(Y\times_{i(X)}1)$.
The down-then-bottom composition sends a map $Y\to X$ to the family $x:1\to X\mapsto f(Y)\times_{j(X)}1$.
The commutation of the square comes from the commutation of $f$ with finite limits, giving a natural equivalence
$f(Y\times_{i(X)}1) = f(Y)\times_{fi(X)}f(1) = f(Y)\times_{j(X)}1$.
\end{proof}

We define a {\it $\Pi$-pretopos} as a pretopos which is locally cartesian closed.
A morphism of $\Pi$-pretopoi is a morphism of pretopoi which is also a morphism of locally cartesian closed categories.

\begin{theorem}
\label{thm:initial-LCCpretopos}
The category $\pif\cS$ is the initial $\Pi$-pretopos.
\end{theorem}
\begin{proof}
Let $\cE$ be a $\Pi$-pretopos.
Then, $\cE$ is in particular a pretopos and we get a unique pretopos morphism $i:\pif\cS\to \cE$ from \Cref{thm:initial}.
The result will be proved if we show that $i$ is a morphism of locally cartesian closed categories.
For any $X$ in $\pif\cS$, we need to show that the pretopos morphism $i_X:(\pif\cS)\slice X \to \cE \slice {iX}$ preserves exponentials.
By \cref{lem:bcoh-diagrams}, the pretopos morphism $i_X$ is equivalent to $i^X : (\pif\cS)^X \to \cE^X$, and by \cref{thm:CohLCC} it suffices to show $i$ preserves exponentials.
This can be shown with the same strategy as in \Cref{prop:CC}.
Let $X$ and $Y$ be two $\pi$\=/finite spaces and $X\index$ a truncated Kan complexes with finite values and with colimit $X$ (\cref{prop:resolution}).
Then we have $Y^X = \lim_n Y^{X_n}$ in $\pif\cS$
and $(iY)^{iX} = \lim_n (iY)^{X_n}$ in $\cE$.
Since $Y$ is $N$-truncated for some $N$ and $i$ is left-exact, then $iY$ is also $N$-truncated and we can use \cref{lem:final2} to reduce both cosimplicial limits to finite limits.
Then we can use that $i:\pif\cS\to \cE$ preserves finite limits, and therefore sends $Y^X = \lim_n Y^{X_n}$ to $\lim_n (iY)^{X_n} = (iY)^{iX}$.
\end{proof}

\begin{remark}
\label{rem:strength}
To appreciate the strength of the initiality condition of \cref{thm:initial-LCCpretopos}, it is useful to compare it with the initial property of the topos $\cS$.
Recall that the category of spaces $\cS$ is initial in the category of topoi and cocontinuous and left-exact functors \cite[Proposition 6.3.4.1]{Lurie:HTT}.
However, $\cS$ is no longer initial in the (non-full) subcategory of topoi and cocontinuous and left-exact functors which are also morphisms of locally cartesian closed categories.
If this was true, this would imply that for any topos $\cE$, the canonical cocontinuous and left-exact functor $i:\cS\to \cE$ always preserves exponentials ($i(Y^X) = (iY)^{iX}$) which is false if $\cE$ is not locally contractible (for example, if $\cE$ is the category of sheaves over the Baire space $\NN^\NN$, the endomorphism sheaf of the constant sheaf $\NN$ is not constant, since the canonical section $ev:\NN^\NN\times \NN\to \NN$ is not locally constant).
However, it is always true that the restriction $\pif\cS \subto \cS\xto i \cE$ does preserve exponentials.
In fact, any topos $\cE$ being a $\Pi$-pretopos, 
\cref{thm:initial-LCCpretopos} says that the canonical functor $\pif\cS \to \cE$ is even a morphism of locally cartesian closed categories.
\end{remark}

\subsection{Initial Boolean \texorpdfstring{$\Pi\Omega$}{Pi-Omega}-pretopos}

The category of finite sets is known to be the universal Boolean elementary 1-topos \cite[pp. 71--73]{Awodey:thesis}.
We can deduce from \cref{thm:initial-LCCpretopos} a similar result for $\pif\cS$.
We define a {\it $\Pi\Omega$-pretopos} as a $\Pi$-pretopos which admits a subobject classifier.
A morphism of $\Pi\Omega$-pretopoi is a morphism of $\Pi$-pretopoi which preserves the subobject classifier.
A $\Pi\Omega$-pretopos is said to be Boolean if its subobject classifier is isomorphic to $2=1+1$.
The category of Boolean $\Pi\Omega$-pretopos is defined as a full subcategory of that of $\Pi\Omega$-pretopoi.

\begin{corollary}
\label{thm:initial-boolean}
The category $\pif\cS$ is the initial Boolean $\Pi\Omega$-pretopos.
\end{corollary}
\begin{proof}
Let $\cE$ be a Boolean $\Pi\Omega$-pretopos.
It is sufficient to prove that the morphism $i:\pif\cS \to \cE$ of $\Pi$-pretopoi given by \cref{thm:initial-LCCpretopos} is in fact a morphism of Boolean $\Pi\Omega$-pretopoi, that is that $i$ preserves the subobject classifiers ($i(\Omega_{\pif\cS}) = \Omega_\cE$).
By assumption the subobject classifier of $\cE$ is $\Omega_\cE=2=1+1$.
Using the fact that $i$ preserves sums and \cref{prop:CohOmega}, we get that $\Omega_\cE = 2 = i(2) = i(\Omega_{\pif\cS})$.
\end{proof}

\begin{appendices}
\section{Univalent families in pretopoi}
\label[appendix]{app:univalence}

This appendix studies univalent families/maps in the context of pretopoi.
We study the notion of univalent reflection of a map and the condition of having enough univalent maps.
The main result is a criteria for the existence of the univalent reflection and a recipe for its construction (\cref{app:main}).
This allows to build all univalent maps.

The construction of the univalent reflection of a map was first done in \cite[Theorem 5.1]{vdBerg-Moerdijk:univalence} (see also \cite[Proposition 5.14]{Stenzel:completion}) but in the setting of model categories.
A treatment of univalence in the \oo categorical setting is done in \cite{Gepner-Kock:univalence,Rasekh:univalence} and most of the material of this appendix appears already in those references but with some differences.
In the first reference, the setting is that of presentable locally cartesian closed categories, and the problems of univalent reflection is not considered.
In the second reference, the setting is that of categories with finite limits only.
The univalent reflection of a map is essentially in Theorem 4.4, but the focus is instead on the construction of the corresponding ``categorical subuniverse'' (see \cref{app:rem:reflection=completion:2}).

\subsection{Pretopoi and their universes}
\label[appendix]{app:sec:pretopos}

We start this appendix with some recollections on descent and a characterization of pretopoi in terms of their universe (\cref{prop:descent-pretopos}).
This will be useful to work with local classes of maps.

\medskip
In category $\cE$ with pullbacks, a class $\cD$ of diagrams (possibly indexed by different categories) is said to be {\it cartesian} if, for every $X\index$ in $\cD$ and every cartesian natural transformation $Y\index\to X\index$, $Y\index$ is also in $\cD$.
Examples of such cartesian class include $I$-indexed diagrams for a fixed collection of categories $I$, but also the class of Segal groupoids (see below).

We fix a class a cartesian class of diagrams $\cD$ such that their colimits exists in $\cE$.
Let $c:X\index \to X$ be a colimit cone where $X\index$ is in $\cD$.
The pullback along $c$ induces a functor $c^*:\cE\slice{|X\index|}\to \lim \cE\slice {X\index}$.
Recall that the category $\lim \cE\slice {X\index}$ is equivalent to that of cartesian natural transformations $Y\index\to X\index$ \cite[Corollary 3.3.3.2]{Lurie:HTT}.
By assumption on $\cD$, every such $Y\index$ is in $\cD$ and has a colimit in $\cE$.
This colimit functor defines the left adjoint $c_! \dashv c^*$.

\begin{definition}[Descent]
\begin{enumerate}
\item The colimit cone $c:X\index\to X$ is called {\it stable} if, for every $Y\to X$, the cone $Y\times_XX\index \to Y$ is a colimit cone.
Equivalently, $c$ is stable 
if and only if $c^*$ is fully faithful,
if and only if the counit of the adjunction $c_! \dashv c^*$ is invertible.

\item The colimit cone $c:X\index\to X$ is called {\it efficient} if, for every colimit cone $Y\index\to Y$ and every cartesian natural transformation $Y\index\to X\index$, the natural transformation $Y\index \to Y\times_XX\index$ is an isomorphism.
Equivalently, $c$ is efficient 
if and only if $c_!$ is fully faithful,
if and only the unit of the adjunction $c_! \dashv c^*$ is invertible.

\item The colimit cone $c:X\index\to X$ has {\it descent} if it is stable and efficient.
Equivalently, $c$ has descent if and only the adjunction $c_! \dashv c^*$ is an equivalence.
\end{enumerate}
\end{definition}

We now specialize the previous results to Segal groupoids.
We fix a category $\cE$ with finite limits.
We shall say that a simplicial diagram $X\index:\Delta\op\to \cE$ is a {\it Segal groupoid} if it is a groupoid object in the sense of \cite[Definition 6.1.2.7]{Lurie:HTT}.
The existence of finite limits in $\cE$ simplifies the definition (see \cref{rem:Lurie-groupoid}): the functor $X[\Lambda^k[n]]$ becomes representable by an object $X^{\Lambda^k[n]}$ in $\cE$ and Condition~(3) of \cite[Proposition 6.1.2.6]{Lurie:HTT} is equivalent to the canonical map $X_n\to X^{\Lambda^k[n]}$ being invertible.

Equivalently, a simplicial diagram $X\index:\Delta\op\to \cE$ is a Segal groupoid if the associated functor $\Map-{X\index}:\cE\op\to \cS^{\Delta\op}$ takes values in the subcategory of Segal groupoids in the sense of \cref{sec:groupoids}.
Segal groupoids verify in particular the Segal conditions $X_n= X_1\times_{X_0}\dots\times_{X_0}X_1$, and we leave the reader the proof that that they form a cartesian class of diagrams.
Let $f:X\to Y$ be a map in $\cE$, and let $X\index$ be its nerve $N(f)\index$ ($X_n = X\times_Y \dots \times_YX$).
Then, $X\index$ is a Segal groupoid.

Let $X\index$ be a simplicial object in $\cE$ with a colimit $|X\index|$.
We shall call $|X\index|$ the {\it quotient} of $X\index$ and refer to the canonical map $q:X_0\to |X\index|$ as the {\it quotient map}.
The quotient of a Segal groupoid is {\it effective} if the canonical map $X\index \to N(q)$ (where $N(q)$ is the nerve of the quotient map) is invertible in $\cE\splx$.
Because of the Segal relations, this condition is equivalent to the following square being a pullback:
\[
\begin{tikzcd}
X_1 \ar[r,"d_0"]\ar[d,"d_1"] & X_0\ar[d]\\
X_0 \ar[r] & {|X\index|}\,.
\end{tikzcd}
\]

\begin{proposition}
\label{prop:effective-efficient}
If quotients of Segal groupoids are stable, 
then a Segal groupoid has an effective quotient if and only if its quotient map is efficient.
\end{proposition}

The proof will need a few lemmas.
We denote by $\mathsf{dec}:\Delta\to \Delta$ the morphism sending $[n]$ to $[n]\star [0]=[n+1]$ where $[n]\star [m] = [n+m+1]$ is the ordinal sum.
The canonical inclusion $[n]\subset [n]\star [0]$ induces a natural transformation $\iota:1\to \mathsf{dec}$ from the identity functor.
The {\it decalage} of simplicial object $X\index:\Delta\op\to \cE$ is the composition 
$X\indexplus = X\index \circ \mathsf{dec}\op :\Delta\op\to\Delta\op\to \cE$.
The transformation $\iota$ induces a natural transformation $X\indexplus\to X\index$.

\begin{lemma}[{\cite[Lemma 6.1.3.17]{Lurie:HTT}}]
\label{lem:split}
If $X\index$ is a Segal groupoid, then 
the natural transformation $X\indexplus\to X\index$ is cartesian and 
$X\indexplus$ is a Segal groupoid whose colimit is $X_0$.
\end{lemma}

We shall say that the decalage of $X\index$ is {\it effective} if the morphism of cocones
\begin{equation}
\label{square:eff-decalage}
\begin{tikzcd}
X\indexplus \ar[r]\ar[d] & X_0\ar[d]\\
X\index \ar[r] & {|X\index|}
\end{tikzcd}
\end{equation}
is a cartesian (in the sense that it is cartesian in the category of simplicial diagrams, when the right hand side is viewed as a constant diagram).

\begin{lemma}
\label{lem:eff-dec}
The quotient of a Segal groupoid $X\index$ is effective if and only if its decalage is effective.
\end{lemma}
\begin{proof}
Since $X\indexplus\to X\index$ is a cartesian natural transformation, the square \eqref{square:eff-decalage} is cartesian if and only if the ``initial'' square
\[
\begin{tikzcd}
X_1 \ar[r]\ar[d] & X_0\ar[d]\\
X_0 \ar[r] & {|X\index|}
\end{tikzcd}
\]
is cartesian.
\end{proof}
\begin{proof}[Proof of \cref{prop:effective-efficient}]

Sufficient condition.
We fix a Segal groupoid $X\index$.
Since the morphism $X\indexplus\to X\index$ is cartesian by \cref{lem:split}, 
the square \eqref{square:eff-decalage} is cartesian morphism of cocones when the quotient map of $X\index$ is efficient.
Then the conclusion follows from \cref{lem:eff-dec}.

Necessary condition.
Using the special shape of simplicial diagrams, 
the quotient map of a Segal groupoid $X\index$ is efficient 
if and only if, for every cartesian natural transformation $Y\index\to X\index$, the following square is a cartesian morphism of cocones
\begin{equation}
\label{eq:effectivity-to-efficiency:1}
\begin{tikzcd}
Y\index \ar[r]\ar[d] & {|Y\index|}\ar[d] \\
X\index \ar[r] & {|X\index|}\,,
\end{tikzcd}
\end{equation}
if and only if, for every cartesian natural transformation $Y\index\to X\index$, the following square is a cartesian in $\cE$
\begin{equation}
\label{eq:effectivity-to-efficiency:2}
\begin{tikzcd}
Y_0 \ar[r]\ar[d] & {|Y\index|}\ar[d] \\
X_0 \ar[r] & {|X\index|}\,.
\end{tikzcd}
\end{equation}

We consider the pullback of the square \eqref{eq:effectivity-to-efficiency:1} along the map $X_0\to |X\index|$
\[
\begin{tikzcd}[sep=small]
&Y\index \ar[rr]\ar[dd] && {|Y\index|}\ar[dd]\\
Y\index\times_{|X\index|}X_0
\ar[rr,crossing over]\ar[dd]  \ar[ur]
&
&
{|Y\index|}\times_{|X\index|}X_0 \ar[ur]
\\
& X\index \ar[rr] && {|X\index|}
\\
X\index\times_{|X\index|}X_0
\ar[rr]  \ar[ur]
&&
X_0 \ar[from=uu,crossing over] \ar[ur]
\,,
\end{tikzcd}
\]
By \cref{lem:eff-dec}, we know that $X\index\times_{|X\index|}X_0 = X\indexplus$.
Using the fact that the natural transformation $Y\index\to X\index$ is cartesian we can compute that $Y\index\times_{|X\index|}X_0 = Y\indexplus$.
Then, by stability of quotients of Segal groupoids, the object ${|Y\index|}\times_{|X\index|}X_0$ is the colimit of the diagram $Y\indexplus$, that is $Y_0$.
Therefore, the right face is the square \eqref{eq:effectivity-to-efficiency:2}.
This shows that it is a cartesian square and that $X\index$ has an efficient quotient map.
\end{proof}

The following definition is a mere reformulation of \cite[Definition A.6.1.1]{Lurie:SAG} using \cref{prop:effective-efficient}.

\begin{definition}[Pretopos]
\label{def:pretopos}
A {\it pretopos} is a category with finite limits, finite sums, and quotients of Segal groupoids such that
\begin{enumerate}[label=\roman*)]
\item sums have descent ($\cE\slice{\coprod X_i}=\prod\cE\slice{X_i}$),
\item Segal groupoids have descent ($\cE\slice{|X\index|}=\lim\cE\slice{X\index}$).
\end{enumerate}
\end{definition}

Consider the two functors $Q:\cE^{\Delta\op} \rightleftarrows \Arr\cE: N$ where $Q$ sends a diagram $X\index$ to the quotient map $q:X_0\to |X\index|$, and $N$ sends a map $f$ to its nerve $N(f)$.
We leave to the reader the proof that $Q$ is left adjoint to $N$.

\begin{lemma}
\label{lem:descent-pretopos}	
The adjunction $Q \dashv N$ restricts to an equivalence between the subcategory of Segal groupoids and the subcategory of surjections.
\end{lemma}
\begin{proof}
Any adjunction restricts to an equivalence between the subcategories of objects for which the unit or counit is invertible.
For $X\index$ a simplicial object the map $X\index \to NQ(X\index)$ is an isomorphism if and only if $X\index$ is a Segal groupoid.
The condition is sufficient since $N(f)$ is a Segal groupoid for every map $f$.
And it is necessary by effectivity of Segal groupoids in $\cE$.
For a map $f:A\to B$, the map $QN(f):A\to C$ is the surjection part of $f$ and $C\to B$ is the image of $f$ (see the proof of \cite[Proposition A.6.2.1]{Lurie:SAG}).
Therefore $QN(f)\to f$ is an isomorphism if and only if 	$f$ is a surjection.
\end{proof}

\begin{definition}[Coherent topology and sheaves]
\label{def:coh-topo}
A finite family $X_i\to X$ of maps is said to be {\it covering} if $\coprod X_i \to X$ is a surjection.
The {\it coherent topology} on $\cE$ (called the effective epimorphism topology in \cite[Corollary A.6.2.3]{Lurie:SAG}) is the Grothendieck topology generated by the finite covering families.

Recall from \cite[Proposition A.6.2.5]{Lurie:SAG} that a functor $F:\cE\op\to \cC$ is a sheaf for the coherent topology if $F(\coprod X_i) = \prod F(X_i)$ for every finite set of objects $X_i$ and if $F(X) = \lim_{\Delta}F(N(f)_n)$ for every surjection $f:Y\to X$.
The equivalence of \cref{lem:descent-pretopos}, show that this last condition is equivalent to asking that $F(|X\index|) = \lim_{\Delta}F(X\index)$ for every Segal groupoid object.
\end{definition}

\medskip

Recall that in a category with pullbacks $\cE$, any map $Y\to X$ can be thought of as a family of objects in $\cE$ indexed by $B$ (which morally are the fibers of the map).
The slice category $\cE\slice B$ can then be thought of as the category of families of objects indexed by $X$, and the corresponding internal groupoid (core) $\cE\slice X\int$ is the groupoid of families of objects indexed by $X$.
We consider the following two functors, classifying families of objects:
\begin{align*}
\sfE_\mathsf{cat}:\cE\op &\tto \Cat \\
X&\mto \cE\slice X\,,\\
\sfE:\cE\op &\tto \cS \\
X&\mto \cE\slice X\int\,.
\end{align*}
The first one is arguably more natural since it contains the whole structure of families of objects, but its values in $\Cat$ makes it slightly harder to deal with.
Restricting the values to groupoids makes it easier to compare it to functor of points of objects of $\cE$.

\begin{definition}[Universe]
\label{def:universe}
Let $\cE$ be a category with pullbacks.
The functor $\sfE_\mathsf{cat}$ is called the {\it categorical universe} of $\cE$,
and the functor $\sfE$ is called the {\it groupoidal universe} of $\cE$.
\end{definition}

\begin{remark}
\label{rem:unstraightening}
\label{rem:univ-fam}
In the previous definition, we favored the functorial point of view over the fibrational one.
The two are of course equivalent using straightening/unstraightening of \cite[\S3]{Lurie:HTT}.
The unstraightening of the categorical universe is the codomain fibration $cod:\Arr\cE\to \cE$ \cite[Lemma~6.1.1.1]{Lurie:HTT}.
In particular, the category of elements of $\sfE_\mathsf{cat}$ is $\Arr\cE$, the arrow category of $\cE$.
Similarly, the fibration corresponding to the groupoidal universe is $cod:\Arr\cE\cart\to \cE$, and the category $\cE\slice\sfE$ of elements of $\sfE$ is equivalent to the wide subcategory $\Arr\cE\cart\subset\Arr\cE$ consisting of only the cartesian morphisms of the codomain fibration (which, here, are simply the cartesian squares).

Let us make the equivalence $\cE\slice\sfE=\Arr\cE\cart$ more precise.
The inclusion $\cE\subset\Sh\cE$ induces a functor $\Arr\cE\cart\to \Arr{\Sh\cE}$, sending a a map $A\to B$ in $\cE$ to the same map viewed in $\Sh\cE$.
The two diagrams obtained by extracting the domain and codomains are related by a cartesian natural transformation.
By straightening, the colimit of the codomain diagram in $\Sh\cE$ is $\sfE$.
We denote $\sfE'$ the colimit of the domain diagram, and $u_\sfE:\sfE'\to \sfE$ the colimit of the whole diagram in $\Arr{\Sh\cE}$.
By descent in $\Sh\cE$, for every map $f:A\to B$ in $\cE$, we have a unique cartesian square 
\[
\begin{tikzcd}
A\ar[r]\ar[d,"f"']\pbmark & \sfE'\ar[d,"u_\sfE"]\\	
B\ar[r] &\sfE\,.
\end{tikzcd}
\]
The equivalence between $\Arr\cE\cart=\cE\slice \sfE$ can now be made more precise: in one direction, it is given by the functor $\Arr\cE\cart\to \cE\slice\sfE$ sending $A\to B$ to the classifying map $B\to \sfE$, and the inverse equivalence $\cE\slice\sfE\to \Arr\cE\cart$ is given by the previous pullback.
\end{remark}

\begin{proposition}
\label{prop:descent-pretopos}
Let $\cE$ be a category with finite limits, finite sums and quotients of Segal groupoids.
The following conditions are equivalent:
\begin{enumerate}
\item\label{prop:descent-pretopos:1} $\cE$ is a pretopos,
\item\label{prop:descent-pretopos:2} the categorical universe $\sfE_\mathsf{cat}$ is a sheaf for the coherent topology,
\item\label{prop:descent-pretopos:3} sums and quotients of Segal groupoids are stable and the groupoidal universe $\sfE$ is a sheaf for the coherent topology.
\end{enumerate}

\end{proposition}
\begin{proof}
\eqref{prop:descent-pretopos:1}$\Leftrightarrow$\eqref{prop:descent-pretopos:2} follows from the description of sheaves.
For \eqref{prop:descent-pretopos:2}$\Leftrightarrow$\eqref{prop:descent-pretopos:3}, we note first that the stability condition says that the functors $c_*$ are fully faithful and the sheaf condition says that $c_!\dashv c^*$ is an equivalence on internal groupoids.
Then, the equivalence follows from the fact that an adjunction $L:C\rightleftarrows D:R$ is an equivalence if and only if $R$ is fully faithful (i.e. if $L$ is a localization) and $L$ and $R$ restricts to an equivalence on the internal groupoids (which implies that the localization $L$ is conservative, thus an equivalence).
\end{proof}

\begin{definition}[Universal family]
\label{app:def:univ-fam}
We shall say that the map $u_\sfE:\sfE'\to \sfE$ is the {\it universal family} of $\cE$.
\end{definition}

\subsection{Local classes of maps and univalent families in pretopoi}
\label[appendix]{app:sec:univ-maps}

Let $\cE$ be a (small) pretopos.
An {\it internal family of objects} in $\cE$ (a {\it family} for short) is simply a map $f:X'\to X$ in $\cE$. 
A morphism $f\to g$ of internal families is a cartesian square in $\cE$
\begin{equation}
\label{eq:cart-mor}
\begin{tikzcd}
X'\ar[d,"f"'] \ar[r] \pbmark & Y'\ar[d,"g"]	\\
X \ar[r] & Y	\,.
\end{tikzcd}
\end{equation}
We shall call such a square a {\it cartesian morphism} of families.
The {\it category of internal families} in $\cE$ is then the category $\Arr\cE\cart=\cE\slice\sfE$ of elements of $\sfE$ (\cref{rem:unstraightening}).

\begin{remark}
The category of families $\Arr\cE\cart=\cE\slice\sfE$ has pullbacks but no terminal object nor finite products.
For every map $f:X'\to X$, the slice $(\Arr\cE\cart)\slice f$ is equivalent to $\cE\slice X$, which is a pretopos.
This shows that $\Arr\cE\cart$ is a {\it local pretopos} in the sense of \cite[Definition A.6.1.1]{Lurie:SAG}.
\end{remark}

\begin{definition}[Cartesian surjection of families]
\label{def:cart-surj}
A cartesian morphism $f\to g$ is a {\it cartesian surjection} if the map $X\to Y$ (and therefore the map $X'\to Y'$) in \eqref{eq:cart-mor} is a surjection.
We denote them by $f\surj g$ for short.
\end{definition}

Cartesian morphisms are called BM-equivalences in \cite[5.2]{Stenzel:completion}.

\medskip
The following lemma is left to the reader.

\begin{lemma}
\label{lem:contraction}
Cartesian surjection are closed under composition, finite sums, and base change in $\Arr\cE\cart$.
\end{lemma}

\begin{definition}[Local classes]
\label{def:local-class}
A class of maps $\cU$ in a pretopos $\cE$ is called {\it local} if
it is closed under base change,
closed under finite sums in the arrow category $\Arr\cE$,
and closed under descent along cartesian surjections: if $f\surj g$ is a cartesian surjection and $f$ in $\cU$, then $g$ is in $\cU$.
We shall often consider local classes as full subcategories of $\Arr\cE\cart$.
We denote by $\Loc\cE$ the poset of local classes ordered by inclusion.
\end{definition}

\begin{remark}
\label{rem:local-def}
This definition is an adaptation of \cite[Proposition 6.2.3.14]{Lurie:HTT} to the context of pretopoi (where only finite coproducts exists).
\end{remark}

\begin{lemma}
\label{lem:local-class=subsheaf}
Local classes of maps are in bijection with subsheaves of $\sfE:\cE\op\to \cS$ for the coherent topology.
\end{lemma}
\begin{proof}
We leave to the reader the proof that subpresheaves of $\sfE$ are in bijection with classes of maps closed under base change.
Given such a class $\cU$, we define a subfunctor $\sfU\subset \sfE$ such that $\sfU(X)\subset \cE\slice X\int$ is the subspace spanned by maps $X'\to X$ that are in $\cU$.
The stability by base change of $\cU$ ensures that this is indeed a subfunctor. 
If $\cU$ is local, the stability by sums and the descent along cartesian surjections are equivalent to the sheaf conditions.
\end{proof}

Recall from \cite[Proposition 6.1.2]{Cisinski}, that a presheaf $\sfU:\cE\op\to \cS$ (or the corresponding right fibration $\cU\to \cE$) is representable (by an object of $\cE$) if and only if its category of elements $\cE\slice\sfU$ has a terminal object $(x,u)$ (the element $u$ provides an isomorphism between the functors $\Map-x\to \sfU$).
We shall say that $\sfU$ is represented by $(x,u)$ or sometimes simply by $x$.

\begin{lemma}
\label{lem:repres-local-class}
The subsheaf $\sfU\subset\sfE$ associated to a local class $\cU$ is representable if and only if the local class has a terminal object when viewed as a full subcategory $\cU\subset \Arr\cE\cart$.
\end{lemma}
\begin{proof}
The full subcategory $\cU\subset\Arr\cE\cart$ is the category of elements of $\sfU$.
The result follows from the fact that a presheaf (or the corresponding right fibration) is representable if and only if its category of elements has a terminal object \cite[Proposition 6.1.2]{Cisinski}.
\end{proof}

Given a map $f$ in $\cE$, we denote by $\{f\}\loc$ the class of maps $g$ for which there exists a span $f\ot h\surj g$ in $\Arr\cE\cart$ where the right leg is a cartesian surjection.

\begin{lemma}
\label{lem:gen-local}
The class $\{f\}\loc$ is the smallest local class containing $f$.
\end{lemma}
\begin{proof}
Let $\cU$ be a local class containing $f$.
Since $\cU$ is closed under base change and descent along cartesian surjections by definition, it must contain $\{f\}\loc$.
Let us see now that $\{f\}\loc$ is local.
We consider a diagram $f\ot h\surj g$.
Let $g'\to g$ be a cartesian morphism.
By \cref{lem:contraction}, the diagram $f\ot h\times_gg'\surj g'$ shows that $g'$ is in $\{f\}\loc$, and that $\{f\}\loc$ is closed under base change.
Let $g\surj g'$ be a cartesian surjection, 
By \cref{lem:contraction}, the diagram $f\ot h\surj g\surj g'$ shows that $g'$ is in $\{f\}\loc$, and that $\{f\}\loc$ is closed under descent along cartesian surjections.
Let $g_i$ be a family of maps in $\{f\}\loc$ indexed by a finite set $I$, and let $f\ot h_i\surj g_i$ the associated diagrams.
By \cref{lem:contraction} again, the diagram $I\times f = \coprod f\ot \coprod h_i\surj \coprod g_i$ (where the first sum is that of the constant diagram with value $f$) is a diagram of cartesian morphisms.
The codiagonal morphism $I\times f =\coprod f\to f$ is a cartesian morphism since it is a pullback along the map $I\to 1$.
The cartesian span $f \ot \coprod f\ot \coprod h_i\surj \coprod g_i$ shows that $\coprod g_i$ is in $\{f\}\loc$, and that $\{f\}\loc$ is closed under finite sums.
This proves that $\{f\}\loc$ is a local class.
\end{proof}

\begin{remark}
\label{rem:infinite-sums}
If $\cE$ has sums indexed by arbitrary sets, and if surjections in $\cE$ are closed under such sums (which is the case if $\cE$ is a topos), the proof shows that the class $\{f\}\loc$ would also be closed under these sums (and be local in the sense of \cite[Proposition 6.2.3.14]{Lurie:HTT}).
\end{remark}

\begin{definition}[Binded local classes]
A local class is {\it binded} if it is of the form $\{f\}\loc$ for some map $f$, such an $f$ is called a {\it bind} for the class.
\end{definition}

\begin{remark}
Every local class in $\cE$ is a suprema of $\{f\}\loc$ in $\Loc\cE$.
Because $\cE$ has finite sums, this suprema is actually a directed union.
The binded local classes are the compact objects of the poset $\Loc\cE$.
\end{remark}

Any cartesian morphism $f\to g$ induces an inclusion $\{f\}\loc \subset \{g\}\loc$.
This defines a functor $L:\Arr\cE\cart \to \Loc\cE$, whose image is the subposet $\BLoc\cE$ of binded local classes.
A cartesian morphism $f\to g$ is called a {\it family equivalence} if $\{f\}\loc = \{g\}\loc$.
This relation is equivalent to the existence of two spans $f\ot h\surj g$ and $g\ot k\surj f$.
By considering $h+k$, we can replace these spans by the single span $f\jrus h+k \surj g$, where the two legs are cartesian surjections.
This description of family equivalences proves the following result.

\begin{lemma}
\label{lem:fam-equiv}
Every cartesian surjection is a family equivalence, and the localization of $\Arr\cE\cart$ inverting family equivalences is generated by inverting cartesian surjections.
\end{lemma}

\begin{lemma}
\label{lem:terminal=bind}
If $u$ is terminal in a local class $\cU$, then it is a bind for $\cU$.
\end{lemma}
\begin{proof}
We have always $\{u\}\loc\subseteq\cU$.
The hypothesis says that for every $f$ in $\cU$ there exists a unique map $f\to u$.
This proves that $\cU\subseteq\{u\}\loc$.
\end{proof}

\begin{remark}
Terminal objects in a local class $\cU$ are called 
{\it classifying morphisms for $\cU$} in \cite[Definition 6.1.6.1]{Lurie:HTT} and 
{\it universal families for $\cU$} in \cite[3.3]{Gepner-Kock:univalence}.
\end{remark}

\begin{proposition}
\label{prop:equiv-univ}
The following conditions on a map $u$ in a pretopos are equivalent:
\begin{enumerate}
\item\label{prop:equiv-univ:local} $u$ is local with respect to cartesian surjections,
\item\label{prop:equiv-univ:term} $u$ is a terminal object in $\{u\}\loc$,
\item\label{prop:equiv-univ:trunc} $u$ is $(-1)$-truncated in $\Arr\cE\cart = \cE\slice \sfE$.\end{enumerate}
\end{proposition}
\begin{proof}
\eqref{prop:equiv-univ:local}$\Rightarrow$\eqref{prop:equiv-univ:term}.
Let $g$ be in $\{u\}\loc$, there exists a diagram $u\ot f\surj g$.
By assumption on $u$, there exists a unique cartesian morphism $g\to u$ such that the composition $f\surj g\to u$ is the cartesian morphism $f\to u$.
This implies that every map in $\{u\}\loc$ is a base change of $u$.
Now, consider the cartesian surjection $g+g\surj g$.
By locality assumption on $u$, the diagonal $\Map g u \to \Map{g+g}u=\Map g u ^2$ is an isomorphism of spaces.
Since we know that $\Map g u$ is not empty, we must have $\Map g u=1$.
This shows that $u$ is terminal in $\{u\}\loc$.


\smallskip
\noindent\eqref{prop:equiv-univ:term}$\Rightarrow$\eqref{prop:equiv-univ:trunc}.
The map $u$ is $(-1)$-truncated in $\Arr\cE\cart$ if and only if for every $g$ in $\Arr\cE\cart$ $\Map g u$ is either empty or a contractible space.
If $g$ is in $\{u\}\loc$, then we have $\Map g u = 1$ by assumption.
If $g$ is not in $\{u\}\loc$, then we must have $\Map g u = \emptyset$ (otherwise $g$ would be a base of $u$, thus in $\{u\}\loc$).

\smallskip
\noindent\eqref{prop:equiv-univ:trunc}$\Rightarrow$\eqref{prop:equiv-univ:local}.
For a cartesian surjection $f\surj g$, we want to show that $\Map g u \to \Map f u$ is invertible.
Since $u$ is $(-1)$-truncated, $\Map g u$ and $\Map f u$ are either empty or a contractible space.
If $\Map fu$ empty, then $\Map g u$ must be empty also.
If $\Map fu$ is non-empty, then it is enough to show that $\Map g u$ is non-empty.
The cartesian surjection $f\surj g$ is the colimit of its nerve in the local pretopos $\Arr\cE\cart$.
Because $u$ is $(-1)$-truncated, all morphisms $f\times_g \dots \times_g f\xto {p_n} f \to u$ must be equal (where $p_n$ is the projection on the $n$-th factor).
This show that the nerve of $f\surj g$ has a natural transformation to the constant simplicial diagram with values $u$.
Then, by taking the colimit, we get a cartesian morphism $g\to u$.
\end{proof}

\begin{definition}[Univalent family]
We shall say that a map $u$ in $\cE$ is a {\it univalent family} (or {\it univalent} for short) if it satisfies the equivalent conditions of \cref{prop:equiv-univ}.
\end{definition}

\begin{remark}
The characterization of univalent maps as terminal objects is \cite[Proposition 3.8~(4)]{Gepner-Kock:univalence}.
And their characterization as $(-1)$-truncated objects is the definition chosen in \cite[Definition 2.1]{Rasekh:univalence} (this choice has the advantage to work when $\cE$ is only a lex category).
\end{remark}

\begin{remark}
\label{rem:Voevodsky}
The universe $\sfE$ classifies uniquely every family in $\cE$: every map $X'\to X$ corresponds to a unique element $X\to \sfE$ via the equivalence $\Arr\cE\cart=\cE\slice\sfE	$, and the map $X'\to X$ is the pullback of the universal family $u_\sfE:\sfE'\to \sfE$ along $X\to \sfE$ (see \cref{rem:univ-fam}).
If $U'\to U$ is a $(-1)$-truncated object, the families that are pullbacks of $U'\to U$ (i.e classified by $U\subset\sfE$) are also classified uniquely.
This is the meaning intended by Voevodsky when he chose the terminology {\it univalent}.
It was meant to complete the terminology {\it versal}, existing in deformation theory, for a map that classifies everything but in a non-unique way. A map is then {\it uni-versal} if it is versal and univalent.
The {\it univalence axiom} of homotopy type theory is then a condition on a fibration in a certain model 1-category to be sent to a univalent family in the \oo category constructed by localization \cite[\S 3]{Kapulkin-LeFanu-Lumsdaine}.
\end{remark}

The following lemma shows that there is a bijection between univalent maps and representable subsheaves of the (groupoidal) universe $\sfE$.

\begin{lemma}
\label{lem:repres-local-class-univ-bind}
The subsheaf $\sfU\subset\sfE$ associated to a local class $\cU$ is representable by a map $u$ if and only if the map $u$ is a univalent bind for $\cU$.
\end{lemma}
\begin{proof}
By \cref{lem:repres-local-class}, the presheaf $\sfU$ is representable by one of its elements $u$ if and only if the map $u$ is a terminal object of $\cU$ when viewed as a full subcategory of $\Arr\cE\cart$.
Using this and the fact that any terminal object $u$ of $\cU$ is always a bind (\cref{lem:terminal=bind}), the statement is equivalent to prove that, if $u$ is a bind for $\cU$ (i.e. if $\cU=\{u\}\loc$), then it is a terminal object if and only if it is univalent.
But this true by definition of univalent maps.
\end{proof}

\begin{lemma}
\label{app:lem:univ-fam=univ}
Every map $f$ in $\{u\}\loc$ is a colimit of a diagram of maps in $\cE$.
\end{lemma}
\begin{proof}

For a map $f:X\to Y$, an object $B$ in $\cE$, and a map $b:B\to Y$ in $\Sh\cE$, we denote $f_b:X_b\to B$ the pullback of $f$ along $b$.
The $f_b$ define a diagram $\cE\slice Y\to \Arr{\Sh\cE}\cart$, indexed by the category of elements of $Y$.
By descent in $\Sh\cE$, the colimit of this diagram is $f$.
Let us see that this is a diagram of maps in $\cE$ when $f$ is in $\{u\}\loc$.

We need to show that, for any $B$ in $\cE$, any map $g:Z\to B$ in $\{u\}\loc$ is in the subcategory $\cE\subset\Sh\cE$.
We fix a span $u\ot h\surj g$.
Let us see first that we can chose $h$ in $\cE$.
By the argument above, the map $h:F\to G$ is the colimit of the cartesian diagram of maps $h_b$.
And the composition $h_b\to h\to u$ shows that the $h_b$ are pullbacks of $u$ along some element of $B\to \sfE$ and therefore in $\cE$ by definition of $\sfE$.
Let $G_0\surj \cE\slice G$ be a surjective functor from a set $G_0$.
The sum $h':=\coprod_{G_0} h_b$ covers the map $h$ and provides a span $u\ot h'\surj g$.
The codomain of $h'$ is a sum of objects in $\cE$ which covers the object $B$.
Since the objects of $\cE$ are quasi-compact in $\Sh\cE$ (they are even coherent \cite[Proposition A.3.1.3~(1)]{Lurie:SAG}), we can extract a finite covering subfamily and get a cartesian surjection $h'':=\coprod_{G_1} h_b\surj g$ where $G_1\subset G_0$ is a finite subset.
The map $h''$ is now in $\cE\subset\Sh\cE$ and we shall use it instead of $h$ in the span $u\ot h \surj g$.

Now, we need to show that the existence of a cartesian surjection $h\surj g$ from a map in $\cE$ forces $g$ to be in $\cE$.
If $h:F\to G$, the codomain of $h\to g$ is a surjection $G\surj B$ in $\cE$.
Its nerve $N(G\to B)$ is a Segal groupoid in $\cE$.
By universality of colimits, the nerve of the domain surjection $F\surj Z$ is the base change of $N(G\surj B)$ along $h$.
Since the inclusion $\cE\subset\Sh\cE$ is closed under finite limits, $N(F\surj Z)$ is a Segal groupoid in $\cE$.
Thus it quotient map $F\surj Z$ is in $\cE$, and this proves the map $g:Z\to B$ is in $\cE$.

Applied to the maps $f_b$ from the beginning, this shows they are in $\cE\subset \Sh\cE$ and that every map $f$ in $\{u\}\loc$ is the colimit of a diagram of maps $f_b$ in $\cE$.
\end{proof}

\begin{proposition}
\label{app:prop:univ-fam=univ}
For every pretopos $\cE$, the universal family $u_\sfE:\sfE'\to \sfE$ is a univalent map in $\Sh\cE$.
\end{proposition}
\begin{proof}
We put $u=u_\sfE$ for short.
We shall show that $u$ is terminal in $\{u\}\loc\subset\Arr{\Sh\cE}\cart$.
By \cref{app:lem:univ-fam=univ}, every map $f$ in $\{u\}\loc$ is a colimit of a diagram of maps $f_b$ in $\cE$.
Using that maps in $\cE$ have a unique map to $u$ in $\Arr{\Sh\cE}\cart$ (\cref{rem:univ-fam,rem:Voevodsky}), we get $\Map f u = \lim \Map {f_b} u = 1$, showing that $u$ is terminal in $\{u\}\loc$, hence univalent.
\end{proof}

\subsection{Univalent reflection}
\label[appendix]{app:sec:univ-refl}

\begin{definition}[Univalent reflection]
\label{def:enough-univalent-maps}
A {\it univalent reflection} for a map $f$ in a pretopos is a terminal object $u_f$ in $\{f\}\loc$.
A pretopos has {\it enough univalent maps} if every map has a univalent reflection.
\end{definition}

The univalent reflection of a map is called its {\it univalent completion} in 
\cite[\S 5]{vdBerg-Moerdijk:univalence} and \cite[Definition 5.12]{Stenzel:completion}.
Univalent reflections may not exist in a pretopos, but we shall see in \cref{app:main} that they always exists if the pretopos is also locally cartesian closed (notion of $\Pi$-pretopos of \cref{sec:pi-pretopos}).
This will make the connection with the original definition of univalent map by Voevodsky.

\begin{remark}
\label{app:rem:reflection=completion}
The univalent reflection is closely related to the completion of Segal spaces/objects, see \cref{app:rem:reflection=completion:2}.
So much, in fact, that {\it univalent} has become (regrettably) a synonym of {\it complete} in the sense of Rezk (e.g. \cite[Theorem 4.4]{Rasekh:univalence} or \cite[Definition 2.11 and \S 4]{Stenzel:completion}).
\end{remark}

\begin{remark}
The existence of image factorization for maps in $\cE\slice \sfE$, shows that $\cE$ has enough univalent maps if and only if every map $f$ is a base change of a univalent map $u$ (no need for the map $f\to u$ to be a cartesian surjection).
To see this, factor this base change into $f\surj u'\hookrightarrow u$, and remark that $u'$ is another univalent map (since subobjects of $(-1)$-truncated objects are $(-1)$-truncated objects).
\end{remark}

\begin{lemma}
\label{app:lem:univ-reflection}
If $g\to f$ is a cartesian morphism, and if $f$ and $g$ have univalent reflections $u_f$ and $u_g$, there is a monomorphism $u_g\to u_f$ in $\Arr\cE\cart$, and the following canonical square commutes
\[
\begin{tikzcd}
g\ar[r]\ar[d] & f\ar[d]\\
u_g\ar[r] & u_f	\,.
\end{tikzcd}
\]
\end{lemma}
\begin{proof}
The cartesian morphism $g\to f$ induces an inclusion $\{g\}\loc\subset\{f\}\loc$.
Since $u_f$ is terminal in $\{f\}\loc$, there exists a unique morphism $u_g\to u_f$ in $\Arr\cE\cart$.
The terminal nature of $u_f$ also show that the square of the statement commutes.
We are left to show that $u_g\to u_f$ is a monomorphism in $\Arr\cE\cart$.
By \cref{prop:equiv-univ}, both $u_g$ and $u_f$ are $(-1)$-truncated objects.
The result follows from the fact that any morphism between $(-1)$-truncated objects is a monomorphism.
\end{proof}

\begin{lemma}
\label{app:lem:equiv-univ-reflection}
If $f'\surj f$ is a cartesian surjection, the univalent reflection of $f$ and $f'$ are isomorphic.
\end{lemma}
\begin{proof}
If $f'\surj f$ is a cartesian surjection, we have $\{f'\}\loc=\{f\}\loc$ by \cref{lem:fam-equiv}, and the terminal objects are uniquely isomorphic in $\Arr\cE\cart$.
\end{proof}

\begin{proposition}
\label{prop:adjoint-ff}
A pretopos has enough univalent maps if and only if the functor $L:\Arr\cE\cart\to \BLoc\cE$, sending a map $f$ to its local class $\{f\}\loc$, is a reflective localization.
\end{proposition}
\begin{proof}
Recall the following general fact about adjunctions \cite[Proposition 5.2.4.2]{Lurie:HTT}.
A functor $L:C\to D$ has a right adjoint if and only if, for every $d$ in $D$, the comma category $C\slice d = C\times_DD\slice d$ has a terminal object $R(d)$. 
Moreover, this left adjoint is fully faithful (i.e. that $LR=id$) if and only if, for every $d$, $R(d)$ belongs to the fiber $C(d) = C\times_D\{d\}$.

We apply this to the functor $L:\Arr\cE\cart\to \BLoc\cE$.
The fiber $\Arr\cE\cart \times_{\BLoc\cE}\{\cU\}$ of $L$ over a local class $\cU$ is the category of its binds (viewed as a full subcategory of $\Arr\cE\cart$).
The comma category $\Arr\cE\cart \times_{\BLoc\cE}\BLoc\cE\slice \cU$ is the class $\cU$ itself (viewed as a full subcategory of $\Arr\cE\cart$).
So the right adjoint to $L$ exists 
if and only if every binded class admits a terminal object 
if and only if every map has a univalent reflection. 
By \cref{lem:terminal=bind} such a terminal object is a bind and the right adjoint must be fully faithful.
This proves the statement.
\end{proof}

\begin{remark}
As a consequence of \cref{prop:adjoint-ff}, when $\cE$ has enough univalent maps, the reflection $L:\Arr\cE\cart \to \BLoc\cE$ is the localization $\Arr\cE\cart$ along the family equivalences (or the cartesian surjections by \cref{lem:fam-equiv}).
\end{remark}


\begin{proposition}
\label{app:prop:univ=sum-univ}
If $\cE$ has enough univalent maps, then its universal family $\sfE'\to \sfE$ is the colimit in $\Arr{\Sh\cE}\cart$ of all the univalent maps of $\cE$.
In particular, the universe $\sfE$ is the filtered union in $\Sh\cE$ of the codomains of all the univalent maps of $\cE$.
\end{proposition}
\begin{proof}
Univalent maps forms a poset since they are the $(-1)$-truncated objects of $\Arr\cE\cart=\cE\slice\sfE$.
When $\cE$ has enough univalent maps this poset is reflective by \cref{prop:adjoint-ff}.
Since right adjoint are final functors, the universal family $\sfE'\to \sfE$ is the colimit of the smaller diagram $\BLoc\cE\to\cE\slice\sfE=\Arr\cE\cart\to \Arr{\Sh\cE}$.
Since all morphisms between univalent maps are monomorphisms, the colimit of  \cref{app:prop:univ=sum-univ} is actually a union.
Moreover, this union is filtered since the poset $\BLoc\cE$ has finite suprema (the supremum of a finite family $u_i$ of univalent maps is the univalent reflection of $\coprod u_i$).
\end{proof}

\begin{remark}[Construction of the univalent reflection]
\label{rem:trunc}

%
By \cref{prop:equiv-univ}, the problems of 
finding a univalent reflection for a map $f:A\to B$ 
or a $(-1)$-truncation for $f$ in $\Arr\cE\cart$, 
or a $(-1)$-truncation for the classifying map $B\to \sfE$ in $\cE\slice \sfE$ are equivalent.
In $\Sh\cE$, the image $U$ of the map $B\to\sfE$ is constructed as the colimit of the nerve $N(B\to \sfE)$.
The category $\cE$ being a pretopos, the object $U$ will be in $\cE\subset\Sh\cE$ if (and only if) the nerve $N(B\to \sfE)$ is a simplicial diagram of objects in $\cE$.
Since this diagram satisfies the Segal conditions, it is in $\cE$ if and only if the objects $N(B\to \sfE)_0$ and $N(B\to \sfE)_1$ are in $\cE$.
The first condition is trivial since $N(B\to \sfE)_0=B$.
By definition of $\sfE$, the object $N(B\to \sfE)_1 = B\times_\sfE B$ is the object of ``symmetries'' of the family $f:A\to B$ corresponding to $B\to\sfE$, that is 
\[
N(B\to \sfE)_1
\ =\ 
B\times_\sfE B
\ =\ 
\mathrm{Iso}_{B\times B}(p_1^*f,p_2^*f)\,.
\]
This equivalence is proven in \cite[Corollary 3.7]{Gepner-Kock:univalence}.
The Segal groupoid $N(B\to \sfE)$ is also the internal groupoid (core) of the Segal object constructed in \cite[Theorem 4.4~(1)]{Rasekh:univalence}.
There it is shown in {Lemma~2.12} that $\mathrm{Iso}_{B\times B}(p_1^*f,p_2^*f)$ is an object of $\cE$ when $p_1^*f=A\times B\to B\times B$ is an exponentiable map in $\cE$.
When this is the case, the Segal groupoid $N(B\to \sfE)$ corresponds to a Segal groupoid in $\Arr\cE\cart$ (i.e. all commutative squares are cartesian)
\[
\Iso f :=
\begin{tikzcd}
\dots
\ar[r, shift left=2]
\ar[r]
\ar[r, shift right=2]
&
\mathrm{Iso}_{B\times B}(p_1^*f,p_2^*f)\times_B A
\ar[r, shift left]
\ar[r, shift right]
\ar[d]
&
A
\ar[d,"f"]
\\
\dots
\ar[r, shift left=2]
\ar[r]
\ar[r, shift right=2]
&
\mathrm{Iso}_{B\times B}(p_1^*f,p_2^*f)
\ar[r, shift left]
\ar[r, shift right]
&
B
\end{tikzcd}
\]
where the bottom row is $N(B\to \sfE)$ and the top row is the canonical action of $\mathrm{Iso}_{B\times B}(p_1^*f,p_2^*f)$ on the object $A$ (identifying isomorphic fibers of $A\to B$).
The quotient is a map in $\cE$, which is the univalent reflection of $f$.

This recovers the classical definition of univalent maps (see \cite[3.2]{Gepner-Kock:univalence} or \cite[Definition 3.1.3]{Kapulkin-LeFanu-Lumsdaine}): a map $f:A\to B$ is univalent 
if and only if the Segal groupoid $N(B\to \sfE)$ is constant with value $B$,
if and only if the ``reflexivity'' map $B\to\mathrm{Iso}_{B\times B}(p_1^*f,p_2^*f)$ is an isomorphism.

This construction makes also more precise the intuitive meaning of univalent maps:
since the groupoid $\Iso f$ acts on $f$ by identifying all isomorphic fibers, the quotient is ``univalent'' in the sense that it contains only a single copy of each of the fibers of $f$.
\end{remark}

The considerations of \cref{rem:trunc} sketch the proof the following result, which is essentially 
\cite[Theorem 5.1]{vdBerg-Moerdijk:univalence},
\cite[Theorem 4.4]{Rasekh:univalence},
and \cite[Proposition 5.14]{Stenzel:completion},
with a formulation inspired from the second reference.

\begin{proposition}
\label{app:main}
The univalent reflection of a map $f$ exists when $p_1^*f=A\times B\to B\times B$ is an exponentiable maps in $\cE$, and it is given by the quotient of the Segal groupoid $\Iso f$.
\end{proposition}
%

We shall be particularly interested in the following consequence.
\begin{corollary}
\label{app:main:cor}
If $\cE$ is a pretopos which is locally cartesian closed, every map has a univalent reflection.
\end{corollary}

\begin{remark}
\label{app:main:rem}
A particular case of \cref{app:main} is proven in \cite[Lemma 6.2]{Gepner-Kock:univalence} when $f:F\to 1$ is associated to an object $F$.
The local class $\{f\}\loc$ is then the class of ``locally trivial $F$-bundles'', that is the maps whose fibers are all isomorphic to $F$.
Such an $f$ is univalent if and only if the group $\Iso f = \mathrm{Aut}(F)$ is trivial, that is if $F$ is a subterminal object in $\cE$.
Its univalent reflection, that is the ``universal $F$-bundle'', is the quotient of $f$ by the action of $\mathrm{Aut}(F)$, that is the map $F/\mathrm{Aut}(F) \to B\mathrm{Aut}(F)$.
\end{remark}

\begin{remark}
\label{app:rem:reflection=completion:2}
Let us make \cref{app:rem:reflection=completion} more precise and explain the connection of the univalent reflection with the completion of Segal spaces.
Given a map $f:A\to B$, the Segal groupoid $N(B\to \sfE)$ of \cref{rem:trunc} is non-complete Segal object in $\cE$.
Recall that the completion of a Segal groupoid is a constant simplicial diagram with values its quotient $U$. 
This object is a subsheaf $U\subset \sfE$ of the groupoidal universe.
Any such subsheaf can be enhanced into a full subsheaf $\sfU_\mathsf{cat}\subset\sfE_\mathsf{cat}$ of the categorical universe.
The sheaf of categories $\sfU_\mathsf{cat}$ can be described as the completion of the Segal object 
\[
\begin{tikzcd}
\dots
\ar[r, shift left=2]
\ar[r]
\ar[r, shift right=2]
&
\mathrm{Hom}_{B\times B}(p_1^*f,p_2^*f)
\ar[r, shift left]
\ar[r, shift right]
&
B
\end{tikzcd}
\]
where $\mathrm{Hom}_{B\times B}(p_1^*f,p_2^*f)$ is the object of all morphisms (not only isomorphisms) between $p_1^*f$ and $p_2^*f$ in $\cE\slice {B\times B}$.
This object is the one considered \cite[Theorem 4.4~(1)]{Rasekh:univalence}.
\end{remark}

\subsection{Univalent families in topoi}
\label[appendix]{app:sec:univ-topos}

Every topos is a pretopos (see \cite[Theorem 1.6.0.6]{Lurie:HTT} and \cite[Example A.6.1.5]{Lurie:SAG}), which is also locally cartesian closed.
So all the considerations of \cref{app:sec:univ-maps,app:sec:univ-refl,app:sec:univ-refl} apply, in particular \cref{app:main:cor}.
However, it is not immediately obvious that this provides the expected notions or results, because the notion of a local class in a topos does not coincide with that of a local class in a pretopos (see \cref{rem:local-def}).
It would be straightforward to adapt all the material of \cref{app:sec:univ-maps,app:sec:univ-refl,app:sec:univ-refl} to the context of topoi and their local classes, but we shall see that this is not necessary to transpose \cref{app:main}.

To distinguish the two notions of local classes, let us call a {\it completely local class} a local class in the sense of \cite[Definition 6.1.3.8]{Lurie:HTT}.
We denote by $\CLoc\cE$ the poset of completely local classes.
Any local class can be completed into a completely local class (by completing it for sums indexed by infinite sets) and this shows that the inclusion $\CLoc\cE\subset\Loc\cE$ is reflective.
The composition $\Arr\cE\cart\to\Loc\cE\to\CLoc\cE$ defines an operation $f\mapsto \{f\}\cloc$ which is the completely local class generated by $f$.
The poset $\BCLoc\cE$ of binded completely local classes is defined as the image of this functor.
Let us say that in a topos, a {\it univalent reflection} for a map $f$ is a terminal object in $\{f\}\cloc$, and that a topos has {\it enough univalent maps} if every map has a univalent reflection in that sense (i.e. if $\Arr\cE\cart\to\BCLoc\cE$ has a right adjoint).

\begin{lemma}
\label{lem:univ-topos=pretopos}
In a topos, we have always $\{f\}\loc=\{f\}\cloc$.
\end{lemma}
\begin{proof}
This is \cref{rem:infinite-sums}.
\end{proof}

\begin{remark}
This might not imply that every local class is completely local, only the binded ones are.
An infinite suprema of binded classes $\{f_i\}\loc$ might not contains the infinite sum $\coprod f_i$, but the suprema of the completely local classes $\{f_i\}\cloc$ does.
\end{remark}

\begin{proposition}
\label{prop:univ-topos=pretopos}
A topos has enough univalent maps as a topos if and only if it has enough univalent maps as a pretopos.
\end{proposition}
\begin{proof}
\Cref{lem:univ-topos=pretopos} says that $\BLoc\cE=\BCLoc\cE$.
Thus $\Arr\cE\cart\to\BCLoc\cE$ has a right adjoint if and only if $\Arr\cE\cart\to\BLoc\cE$ has a right adjoint.
\end{proof}

\begin{corollary}
\label{cor:univ-topos}
A topos has always enough univalent maps.
\end{corollary}
\begin{proof}
A topos is always a locally cartesian closed category.
By \cref{app:main:cor} it has enough univalent maps as a pretopos.
Then the result follows from \cref{prop:univ-topos=pretopos}.
\end{proof}

\subsubsection{Univalent families in $\cS$}
\label{app:sec:univ-fam-S}

\Cref{cor:univ-topos} applies in particular to $\cS$.
Moreover in this example, one can take advantage of the fact that every space $B$ can be covered by a set $B_0$ to give a description of all univalent maps.
The pullback of a map $f:A\to B$ along a cover by a set $B_0\to B$ is a map $f'$ which is a sum of maps $F\to 1$ for some spaces $F$ (which are fibers of the map $f$).
Since $B_0\to B$ is surjective, the univalent reflection of $f$ coincide with that of $f'$.
And, following \cref{app:main:rem}, the local class $\{f\}\loc=\{f'\}\loc$ is that of maps whose fibers are isomorphic to some $F$.
In fact, to ensure that $f'$ and $f$ generate the same local class, and have the same univalent reflection, we do not need the map $B_0\to B$ to be surjective: if $f$ has isomorphic fibers over two connected components of $B$, it is enough for $B_0\to B$ to map into one of them only.

\medskip
Let us detail the case where $f'$ is the sum of two maps $F\to 1$ and $G\to 1$.
If $F$ and $F'$ are not isomorphic, then the quotient of $\Iso {f'}$ is the sum of the univalent reflection of $F\to 1$ and $G\to 1$, that is the sum of $F/\mathrm{Aut}(F) \to B\mathrm{Aut}(F)$ and $G/\mathrm{Aut}(G) \to B\mathrm{Aut}(G)$.
If $F$ and $G$ are isomorphic, then $\Iso {f'}$ identifies them and the quotient is simply $F/\mathrm{Aut}(F) \to B\mathrm{Aut}(F)$.

\medskip
This describes all the univalent maps of $\cS$ as the sums of maps $F/\mathrm{Aut}(F) \to B\mathrm{Aut}(F)$, where $F$ runs in a set of spaces which are {\it two by two non-isomorphic}.
The whole (groupoidal) universe of $\cS$ is then the non-small sum of all $F/\mathrm{Aut}(F) \to B\mathrm{Aut}(F)$ where $F$ runs in a set of representative for all the isomorphism classes of objects in $\cS$.
(Such a description of univalent maps does not work in an arbitrary topos or pretopos since it relies on the ability of covering an arbitrary object by global sections).

\begin{example}[Universe of finite sets]
\label{app:ex:fin-set}
A classical example of a univalent family is the universal family of finite sets. 
We put $\underline n := \{0,\dots,n-1\}$ for $n>0$, and $\underline 0=\emptyset$.
Then, the universal family of finite sets can be defined as the univalent reflection of the map $\coprod_{n\in \NN} \{0,\dots,n-1\}\to \NN$,
or, directly, as the map $u_\set:\coprod_n \underline n/\mathfrak S_n \to \coprod_n B\mathfrak S_n$ (where $\mathfrak S_n$ is the group of permutations of $n$ elements).
\end{example}

\subsection{Univalent families with structure}
\label[appendix]{app:sec:univ-structure}

This section defines the structures of having path identity types, dependent sums, dependent products, and quotients of Segal groupoids, on a local class/subunivers/univalent families.

\subsubsection{Identity types}
\label[appendix]{app:sec:identity-types}

The following definition is taken from \cite[Proposition 15]{Awodey:natural}.
It captures the idea that a family of objects is closed under the construction of path spaces (aka identity types).

\begin{definition}[Family with diagonals]
A local class $\cU$ is said to be {\it closed under diagonals} if, for every $f$ in $\cU$, the map $\Delta f$ is also in $\cU$.
A family $f:A\to B$ is said to be {\it closed under diagonals} if $\Delta f$ is a pullback of $f$.
\end{definition}

\begin{lemma}
\label{lem:pb-diago}
Any cartesian morphism $f'\to f$ induces a cartesian morphism $\Delta f'\to \Delta f$.
Moreover, if $f'\surj f$ is a cartesian surjection, then so is $\Delta f'\to \Delta f$.
\end{lemma}
\begin{proof}
Given a cartesian square
\[
\begin{tikzcd}
A'\ar[r]\ar[d,"f'"']\pbmark &A\ar[d,"f"]\\
B'\ar[r]&B
\end{tikzcd}
\]
The map $\Delta f' \to \Delta f$ is the square
\[
\begin{tikzcd}
A' \ar[r]\ar[d] & A \ar[d] \\    
A' \times_{B'} A' \ar[r] & A \times_BA \,.
\end{tikzcd}
\]
Using that $A'= A\times_BB'$, the previous square becomes
\[
\begin{tikzcd}
A \times_BB' \ar[r]\ar[d] & A \ar[d] \\    
A \times_BA \times_BB' \ar[r] & A \times_BA
\end{tikzcd}
\]
which is clearly cartesian.

If the map $B'\to B$ is a surjection, then so is the map $A \times_BA \times_BB' \to A \times_BA
$ in the previous diagram. 
This proves the second assertion.
\end{proof}

\begin{proposition}
\label{prop:closure-diagonal}
If a family $f$ is closed under diagonals then the local class $\{f\}\loc$ is closed under diagonals.
The converse is true if $f$ is univalent.
\end{proposition}
\begin{proof}
Let $g$ be in $\{f\}\loc$.
By definition of $\{f\}\loc$, there exists a span of cartesian morphisms $g\ot h\surj f$.
By \cref{lem:pb-diago}, we deduce a span $\Delta g\ot \Delta h\surj \Delta f$.
Since $\Delta f$ is in $\{f\}\loc$ and local classes are closed by base change and descent along cartesian surjections, this shows that $\Delta g$ is in $\{f\}\loc$.

For a general $f$, if the class $\{f\}\loc$ is closed under diagonals, we have a span $\Delta f\ot h \surj f$, but we may not have a cartesian morphism $\Delta f \to f$.
However, when $f$ is univalent, such a cartesian morphism always exists since $f$ is terminal in $\{f\}\loc$.
\end{proof}

\begin{lemma}
\label{lem:univ-diagonal}
In a pretopos with enough univalent families, 
if a family $f$ is closed under diagonals, then so is its univalent reflection.
\end{lemma}
\begin{proof}
We fix a cartesian morphism $\Delta f \to f$	, which exists by assumption.
Let $f\to u_f$ be the univalent reflection of $f$.
By \cref{lem:pb-diago}, the cartesian surjection $f\to u_f$ induces a cartesian surjection $\Delta f \to \Delta u_f$.
By \cref{app:lem:univ-reflection} we get a diagram of cartesian morphisms
\[
\begin{tikzcd}
f \ar[d]& \ar[l] \Delta f \ar[r, two heads] \ar[d]& \Delta u_f \ar[d]\\
u_f & \ar[l] u_{\Delta f} \ar[r,"="] & u_{\Delta u_f}\,.
\end{tikzcd}
\]
By \cref{app:lem:equiv-univ-reflection} the map $u_{\Delta f}\to u_{\Delta u_f}$ is an isomorphism.
The resulting map $\Delta u_f\to u_f$ shows that $u_f$ is closed under diagonals.
\end{proof}

\subsubsection{Dependent sums and products}
\label[appendix]{app:sec:dep-sum-products}

We now assume that the pretopos $\cE$ is locally cartesian closed.
If $f:A\to B$ is a map in $\cE$, there exists a triple adjunction $f_!\dashv f^*\dashv f_*$ where $f^*:\cE\slice B\to \cE\slice A$ is the pullback along $f$, and where $f_!$ is given by the composition with $f$.

Let $u:U'\to U$ be a univalent family in $\cE$. 
We consider the associated local class $\cU=\{u\}\loc$, 
and the corresponding subuniverse $\sfU\subset\sfE:\cE\op\to\cS$, 
for which $\sfU(A)\subset (\cE\slice A)\int$ is the subspace spanned by maps $A'\to A$ in $\cU$ (since $u$ is univalent these are the maps that are base change of $u$).
For every map $f:A\to B$ in $\cE$, the base change along $f$ restricts to a functor $f^*:\sfU(B) \to \sfU(A)$.

\begin{definition}[Dependent sums and products]
\label{def:dep-sigma-pi}
Let $u:U'\to U$ be a univalent family in $\cE$.
\begin{enumerate}
\item\label{def:dep-sigma-pi:1} The map $u$ is {\it closed under dependent sums} if, for any map $f:A\to B$ in $\cU$, the left adjoint $f_!$ sends $\sfU(A)$ into $\sfU(B)$.
\item\label{def:dep-sigma-pi:2} The map $u$ is {\it closed under dependent products} if, for any map $f$  in $\cU$, the right adjoint $f_*$ sends $\sfU(A)$ into $\sfU(B)$.
\end{enumerate}
\end{definition}

\begin{remark}
A more general definition of having dependent sums or products for a cartesian fibration $\cF\to \cE$ (or for the corresponding functor $\sfF:\cE\op\to \Cat$) involves Beck--Chevalley conditions \cite{Streicher:fibrations,Anel-Weinberger}, but these conditions are for free here.
Indeed, we are only working with subfibrations of the codomain fibrations where the Beck--Chevalley conditions are always true: for the functors $f_!$, they reduce to the cancellation property of pullback squares, and the condition for $f_*$ follows formally by adjunction.
\end{remark}

\begin{lemma}
\label{lem:sigma=compo}
A univalent map $u:U'\to U$ is closed under dependent sums if and only if the corresponding local class $\cU$ is closed under composition.
\end{lemma}
\begin{proof}
Given a map $f:A\to B$ in $\cU$, the left adjoint $f_!:\cS\slice A\to \cS\slice B$ is the composition by $f$.
It sends $\sfU(A)$ into $\sfU(B)$ if and only if $\cU$ is closed under composition by $f$.
It follows that $f_!$ sends $\sfU(A)$ into $\sfU(B)$ for every $f$ in $\cU$ if and only if $\cU$ is closed under composition.
\end{proof}

\subsubsection{Finite sums and products}
\label[appendix]{app:sec:sum-products}

Let $\cE$ be a pretopos, $\sfE$ its (groupoidal) universe, and $\sfU\subset\sfE$ a subuniverse (a subsheaf for the coherent topology).
For every natural number $n$, we define the {\it universe of families of cardinal $n$} in $\sfU$ as the functor $\sfU^n:\cE\op\to \cS$ where $\sfU^n(A)=\sfU(A)^n$.
When $\sfU=\sfE$, the existence of finite sums and products in $\cE$ defines two natural transformations $+^{(n)},\times^{(n)}:\sfE^n\to \sfE$.

\begin{definition}[Finite sums and products]
\label{def:fin-sigma-pi}
We define only the closure properties for finite sums.
The definitions for finite products are similar and left to the reader.

\begin{enumerate}
\item\label{def:fin-sigma-pi:subuniv} A subuniverse $\sfU\subset\sfE$ is said to be {\it closed under finite sums} if, for every $n$, the functor $+^{(n)}:\sfE^n\to \sfE$ sends $\sfU^n\subset\sfE^n$ to $\sfU\subset\sfE$.

\item\label{def:fin-sigma-pi:local} A local class $\cU$ is said to be {\it closed under finite sums} if, for every finite set of maps $A_i\to B$ in $\cU$, the map $\coprod_i A_i\to B$ is in $\cU$.

\item\label{def:fin-sigma-pi:univ} A univalent map $u:U'\to U$ is said to be {\it closed under finite sums} if for every finite family $A_i\to B$ of pullbacks of $u$ over a common base, the sum map $\coprod_i A_i\to B$ is a pullback of $u$.
\end{enumerate}
\end{definition}

In particular, the closure under empty sum implies that the map $\emptyset\to 1$ is in $\cU$, and by stability by base change, every map $\emptyset \to A$ for $A$ in $\cE$.
Similarly, the closure under empty products implies that the identity map of $1$ is in $\cU$, and, by base change, also every isomorphism of $\cE$.

\begin{lemma}
\label{app:lem:fin-sum-prod}
A subuniverse $\sfU$ is closed under finite sums (finite products)
if and only if the corresponding local class $\cU$ is closed under finite sums (finite products),
if and only if the corresponding univalent map $u$ is closed under finite sums (finite products).
\end{lemma}
\begin{proof}
The equivalences follow from the fact that 
a map $X_i\to A$ is an element of the presheaf $\sfU$ 
if and only if it is in the class $\cU$ (\cref{lem:local-class=subsheaf})
if and only if it is a pullback of $u$ (\cref{prop:equiv-univ}~\eqref{prop:equiv-univ:term}).
\end{proof}

%

\subsubsection{Segal groupoids}
\label[appendix]{app:sec:Segal}

Let $\cE$ be a pretopos, $\sfE$ its (groupoidal) universe, and $\sfU\subset\sfE$ a subuniverse (a subsheaf for the coherent topology).
We define the {\it universe of Segal groupoids} in $\sfU$ as the functor $\sfU^{\Delta\op}:\cE\op\to \cS$ where $\sfU^{\Delta\op}(A)$ is the subgroupoid of $((\cE\slice A)^{\Delta\op})\int$ spanned by those diagrams $X\index$ which are Segal groupoids and such that the maps $X_n\to A$ are elements of $\sfU$.
When $\sfU=\sfE$, the existence of quotients for Segal groupoids defines a natural transformation $q:\sfE^\mathsf{Segal}\to \sfE$.

\begin{definition}[Closure under quotient]
\label{app:def:Segal}
\begin{enumerate}
\item\label{app:def:Segal:subuniv} A subuniverse $\sfU\subset\sfE$ is said to be {\it closed under quotients of Segal groupoids} if the quotient functor $q:\sfE^\mathsf{Segal}\to \sfE$ sends $\sfU^\mathsf{Segal}\subset\sfE^\mathsf{Segal}$ to $\sfU\subset\sfE$.

\item\label{app:def:Segal:local} A local class $\cU$ is said to be {\it closed under quotients of Segal groupoids} if for every Segal groupoid $X\index\to A$ for which the maps $X_n\to A$ are in $\cU$, the quotient map $|X\index|\to A$ is in $\cU$.

\item\label{app:def:Segal:univ} A univalent map $u:U'\to U$ is said to be {\it closed under quotients of Segal groupoids} if for every Segal groupoid $X\index\to A$ for which the maps $X_n\to A$ are pullbacks of $u$, the quotient map $|X\index|\to A$ is a pullback of $u$.

\end{enumerate}
\end{definition}

The unravelling of the first definition shows that the second one is a mere reformulation.
An argument similar to that of \cref{app:lem:fin-sum-prod} proves the following lemma.

\begin{lemma}
\label{app:lem:Segal-gpd}
A subuniverse $\sfU$ is closed under quotients of Segal groupoids 
if and only if the corresponding local class $\cU$ is closed under quotients of Segal groupoids,
if and only if the corresponding univalent map $u$ is closed under quotients of Segal groupoids.
\end{lemma}

\end{appendices}


\providecommand{\bysame}{\leavevmode\hbox to3em{\hrulefill}\thinspace}

\end{document}